\documentclass[final,1p,times]{elsarticle}
\makeatletter
\def\ps@pprintTitle{%
	\let\@oddhead\@empty
	\let\@evenhead\@empty
	\def\@oddfoot{}%
	\let\@evenfoot\@oddfoot}
\makeatother

\usepackage{amssymb}
\usepackage{amsthm}
\usepackage{cleveref}
\usepackage{amsmath} 
\usepackage{mathrsfs}
\usepackage{float}
\usepackage{titlesec}
\usepackage{graphicx}
\usepackage{subfigure}
\usepackage{bm}
\usepackage{color}

\theoremstyle{definition}
\newtheorem{exmp}{Example}[section]

\usepackage[usenames,dvipsnames,svgnames,table]{xcolor}

\newcommand{\bsx}{{\boldsymbol x}} 
\newcommand{\bsy}{{\boldsymbol x'}} 
 
\newcommand{\bsnu}{{\boldsymbol \nu}} 
\newcommand{\bsxi}{{\boldsymbol \xi}} 

\newcommand{\mcA}{{\mathcal A}}
\newcommand{\mcE}{{\mathcal E}}
\newcommand{\mcP}{{\mathcal P}}

\begin{document}

\begin{frontmatter}

\title{Improved convergence of fast integral equation solvers for acoustic scattering by inhomogeneous penetrable media with discontinuous material interface}

\author[label1]{Ambuj Pandey}
\ead{ambujpandey84@gmail.com}
\author[label2]{Akash Anand}
\ead{akasha@iitk.ac.in}
\ead[url]{http://home.iitk.ac.in/~akasha}
\address[label1]{Computing and Mathematical Sciences, Caltech, Pasadena, CA 91125, USA}
\address[label2]{Mathematics and Statistics, IIT Kanpur, Kanpur, UP 208016, India}
\begin{abstract}
In recent years, several fast solvers for the solution of the Lippmann-Schwinger integral equation that mathematically models the scattering of time-harmonic acoustic waves by penetrable inhomogeneous obstacles, have been proposed. While many of these fast methodologies exhibit rapid convergence for smoothly varying scattering configurations, the rate for most of them reduce to either linear or quadratic when material properties are allowed to jump across the interface. A notable exception to this is a recently introduced Nystr\"{o}m scheme [J. Comput. Phys., 311 (2016), 258--274]
that utilizes a specialized quadrature in the boundary region for a high-order treatment of the material interface. In this text, we present a solution framework that relies on the specialized boundary integrator to enhance the convergence rate of other fast, low order methodologies without adding to their computational complexity of $O(N \log N)$ for an $N$-point discretization. In particular, to demonstrate the efficacy of the proposed framework, we explain its implementation to enhance the order to convergence of two schemes, one introduced by Duan and Rokhlin [J. Comput. Phys., 228(6) (2009), 2152--2174]
that  is based on a pre-corrected trapezoidal rule while the other by Bruno and Hyde [J. Comput. Phys., 200(2) (2004), 670--694] which relies on a suitable decomposition of the Green's function via Addition theorem. In addition to a detailed description of these methodologies, we also present a comparative performance study of the improved versions of these two and the Nystr\"{o}m solver in 
[J. Comput. Phys., 311 (2016), 258--274]
through a wide range of numerical experiments.
\end{abstract}


\end{frontmatter}


\section{Introduction}
\label{sec::intro}

Owing to a wide range of applications of acoustic or electromagnetic scattering by penetrable inhomogeneous media, such as, medical imaging, radar, sonar, underwater acoustics \cite{bayliss1985numerical}, quest for a fast and accurate numerical method continues to be of interest, and in fact, several promising methodologies already exist in the literature. Most existing numerical solution schemes broadly fall into two categories, namely, the ones that rely on direct discretization of underlying partial differential equation or their variational counterparts, and the others that work with their integral equation reformulations. 
In particular, the integral equation based approach for the solution of forward scattering problem has been an active area of research that has seen a lot of progress in recent years where several fast and accurate numerical solvers have been proposed. The primary advantage in working with integral equation formulations over its differential equation or variational counterparts \cite{li2010coupling,meddahi2003computing,zwamborn1992three,bayliss1985accuracy,medvinsky2013high} lies in the fact that numerical solution can be made to satisfy the radiation condition simply by suitably employing the radiating fundamental solution. This, thus, avoids complications that otherwise arise where a relatively large computational domain containing the scatterer must be used, together with appropriate absorbing boundary conditions on the boundary of the computational domain \cite{gan1993finite}. 

%

\subsection{Problem statement}

The two dimensional forward scattering problem that we consider in this text is described as follows:
given an obstacle $\Omega$, a bounded open subset of $\mathbb{R}^2$, with a smooth boundary $\partial\Omega$, and an 
incident time-harmonic acoustic wave $u^{inc}$ satisfying
\begin{equation} \label{FHE}
\Delta u^{inc}(\bsx) + \kappa^2 u^{inc}(\bsx) = 0, \ \ \bsx \in \mathbb{R}^2,
\end{equation}
where $\kappa = \omega / c_0$ is the wavenumber, $\omega$ is the angular frequency, and $c_{0}$ is the constant speed of wave outside the inhomogeneity $\Omega$,
find the total acoustic field  
$u$ that satisfies {\cite{colton2012inverse}}
\begin{equation} \label{HE}
\Delta u(\bsx) + \kappa^2 n^2(\bsx) u(\bsx) = 0, \ \  \bsx \in \mathbb{R}^2,
\end{equation}
with the refractive index $n(\bsx) = c_{0} / c(\bsx)$, where $c$, the speed of acoustic wave, is allowed to vary with position within $\Omega$ 
and the scattered field $u^s := u - u^{inc}$ satisfies Sommerfeld radiation condition 
\begin{equation} \label{eq:-Sommerfeld}
\lim_{r \to \infty} \sqrt{r} \left( \frac{\partial u^s}{\partial r} - i \kappa u^s \right) = 0,
\end{equation}
where $r = (x_1^2+x_2^2)^{1/2}$ and $i = \sqrt{-1}$ is the imaginary unit. 

It is well known \cite{colton2012inverse,martin2003acoustic} that the solution $u$ to this scattering problem can be obtained by solving the equivalent 
\textit{Lippmann-Schwinger} integral equation 
given by
 \begin{equation} \label{eq:-Lippmann}
 u(\bsx) + \kappa^2  \mcA(mu)(\bsx) = 
 u^{inc}(\bsx), \hspace{3mm} \hspace{3mm} \bsx \in \mathbb{R}^{2},
\end{equation}
with the volume potential $\mcA$ given by
\begin{equation} \label{eq:-vol-pot}
\mcA(v)(\bsx) = \int \limits_{\Omega}G_{\kappa}(\bsx - \bsy) v(\bsy)d\bsy,
\end{equation}
where
$ 
G_{\kappa}(\bm{x}) = \frac{i}{4}H^{1}_{0}(\kappa|\bsx|),
$
is the radiating fundamental solution of Helmholtz equation in the free space and $m(\bm{x}) = 1-n^{2}(\bsx)$.

\subsection{Overview}
In recent years, a  number of  algorithms, including direct and iterative solvers, have been proposed for the solution of Lippmann-Schwinger equation. While we do not review all such contributions, some recent numerical methods include \cite{chen2002fast,aguilar2004high,bruno2004efficient,duan2009high,hyde2002fast,andersson2005fast,sifuentes2010preconditioned,vainikko2000fast,lanzara2004numerical,anand2007efficient,gillman2014spectrally,ambikasaran2016fast,Anand2015highorder,vico2016fast,egidi2007efficient,amar1983numerical,hesford2010fast,corona2015n,bruno2005higher,bruno2003fast,liu2000high,polimeridis2014stable}. Most fast algorithm, among the cited methods, while converging rapidly for smooth scattering media, yield only linear convergence in the presence of discontinuous scattering media. For instance, schemes introduced in \cite{aguilar2004high,duan2009high}, provide a fast high-order method for smooth scattering media by means of a Fast Fourier Transform (FFT based pre-corrected trapezoidal rule) 
but 
fails to produce high-order accuracy for scattering configurations that contain discontinuous material interfaces. We must note that there do exist fast numerical techniques that fair better in terms of convergence rate while dealing with non-smooth scattering objects, for instance, see \cite{hyde2002fast,andersson2005fast,anand2007efficient,Anand2015highorder}. Among these, the approach in \cite{hyde2002fast,andersson2005fast} converges quadratically in the presence of material discontinuity while the algorithm presented in \cite{anand2007efficient}, though high-order convergent, is computationally well suited only for thin inhomogeneities.  While a more recent contribution, especially designed for a high-order treatment of  discontinuous material interfaces \cite{Anand2015highorder} does converge rapidly with optimal computational cost,
the griding strategy used therein, in certain cases, can restrict the methodology from achieving the theoretical computational complexity. Indeed, as we show later in this text, the scheme in \cite{Anand2015highorder} produces far less accurate approximations to the solution when compared to their counterparts obtained through the proposed approach.  

Apart from these,  a couple of fast direct solvers with computational cost $O(N^{3/2})$ have also been proposed in \cite{gillman2014spectrally,ambikasaran2016fast} that primarily rely on quad-tree data structure. While these methods offer several advantages, such as, robustness and capability to handle some large scale frequency regime, they are designed only for material properties that are globally smooth.

The primary aim of this paper is to provide a framework that allows us to enhance the rate of convergence for those existing $O(N \log N)$ solvers that converge rapidly for smooth inhomogeneities but yield low order in the context of discontinuous scattering media, without  adversely effecting their asymptotic computational complexity. 
Indeed, most fast solvers converge slowly in the presence of material discontinuities because evaluation of the integral operator (\ref{eq:-vol-pot}) through fast algorithms (e.g., FFT, AIM, FMM, etc.) require them to integrate across the material interface. The proposed framework successfully overcomes this difficulty by carving out a certain ``thin" boundary region of the scattering inhomogeneity which, in turn, facilitates decoupling of the problem of integration near material discontinuities from rest of the ``base volume" where variations in material properties are smooth. 
Moreover, this decoupling is affected in such a manner that the base integral has an integrand that
smoothly vanishes at the material boundary thus enabling the aforementioned fast convolution strategies to produce rapidly convergent approximations. To demonstrate this idea, we have implemented two different methods for the integration over base volume; namely, an Addition theorem based algorithm \cite{hyde2002fast} (that we refer to in this text as ATM) and another based on pre-corrected trapezoidal rule method \cite{aguilar2004high,duan2009high} (that we refer to as PTM). We recall that in ATM, high-order convergence is achieved by replacing the contrast function $m(\bsx)$ by its truncated Fourier series while the kernel $G_{\kappa}(\bsx -\bsy)$ in (\ref{eq:-vol-pot}) is replaced by a truncated series given by the Addition theorem in polar coordinates \cite{colton2012inverse}. A high-order accuracy in the  pre-corrected trapezoidal rule based method, on the other hand, is attained by appropriately modifying quadrature weights at a few points in the vicinity of singularity. While the computation of these modified weights require the solution of a severely ill conditioned linear system and is time consuming, the methodology remains an attractive alternative as modified weights need to carefully precomputed and tabulated only once for further use \cite{aguilar2004high,marin2014corrected}.


Rest of the paper is organized as follows:  section \ref{sec:-OS} provides an overview of our algorithm where we explain the use of a  smooth cut-off function to break the evaluation of volume potential (\ref{eq:-vol-pot}) into integrals over boundary and base regions. In this section, we also outline our strategy of breaking the overall computation of integral operator (\ref{eq:-vol-pot}) into four different subproblems, namely {\em base-base}, {\em base-boundary}, {\em boundary-base} and {\em boundary-boundary} interactions.
 In sections \ref{sec:-BubI}-\ref{ref:-bulk-boundary}, we provide a detailed description of aforementioned interactions. A brief account of overall computational complexity of this framework is discussed in section \ref{sec:-Comp-Cost}.
 A variety of numerical results, to validate the performance and accuracy of the two methods implemented under this framework are presented in section \ref{sec:-Num-Res} while conclusions are put forth in section \ref{sec:-Conclu}. 
 

\section{Outline of the scheme}
\label{sec:-OS}
Our numerical approach consists to two basic components, namely,
\begin{enumerate}
\item
approximate operators $\mcA_h(v)$ that converge rapidly to the integral $\mcA(v)$ in (\ref{eq:-vol-pot}) as $h \to 0$ for functions $v$ that are smooth within $\Omega$, where $h$ measures the maximum spacing of the points on the underlying computational grid, say $\Omega_h$, and
\item
the iterative linear algebra solver GMRES \cite{saad1986} for solution of the system of linear equations
\[
u_h(\bsx_\ell) + \kappa^2\mcA_h (mu_h)(\bsx_\ell) = u^{inc}(\bsx_\ell), \ \ \ \bsx_\ell \in \Omega_h,
\] 
resulting from the approximation of equation (\ref{eq:-Lippmann}).
\end{enumerate}
This text focuses on the first component where we introduce an efficient scheme that produces high-order accurate approximations to the volume integral applicable to the Lippmann-Schwinger equation in two dimensions. 

As mentioned in the introduction, the proposed framework achieves high-order convergence by 
suitably breaking the problem of integration in $\mathcal{A}(mu)$ into four subproblems.
Toward this, we being by isolating a thin boundary region of $\Omega$ that we discuss next. 

\subsection{The boundary region}

\begin{figure}[t]
\begin{center}
\subfigure[An illustration of the projection operator $\mcP$]{\includegraphics[scale=0.7]{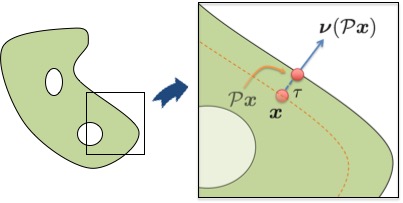} 
\label{fig:-projection}
} 
\subfigure[The cut-off function $\eta_{\tau_0}$]{\includegraphics[scale=0.72]{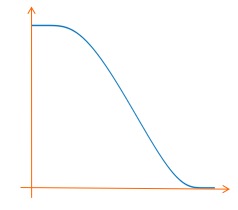} 
\label{fig:-cutoff}
} 
\label{fig:-projection-cutoff}
\end{center}
\caption{Boundary region}
\end{figure}

Since $\partial\Omega$ is assumed to be smooth, there is a neighborhood $N(\partial\Omega)$ such that the orthogonal projection operator $\mcP : N(\partial\Omega) \to \partial\Omega$ satisfying
\[
\mcP\bsx = \bsx + \tau \bsnu(\mcP\bsx),
\] 
where $\bsnu(\bsx)$ denotes the outward unit normal at $\bsx \in \partial\Omega$, is well defined (see Figure \ref{fig:-projection} for an illustration). 
In other words, there exists $\tau_0 > 0$ such that for all $\tau \in [0, \tau_0)$, we can introduce curves $\partial\Omega_{\tau}$ parallel  to $\partial\Omega$ given by the representation
\[
\partial\Omega_{\tau} = \{ \bsx \in \Omega \ :\  (\mcP\bsx - \bsx)\cdot \bsnu(\mcP\bsx) = \tau \}.
\]
Clearly, $\partial\Omega_0 \equiv \partial\Omega$.
The boundary region $\Omega_B$ can thus be defined as the union of these parallel curves, that is,
\[
\Omega_B = \bigcup_{\tau \in (0,\tau_0)} \partial\Omega_{\tau}.
\]

\begin{figure}[!h] 
	\begin{center}	
		\subfigure[A bean shaped scatterer $\Omega$ with a corresonding boundary region $\Omega_{B}$.]	{\includegraphics[clip=true, trim=0 0  00 0, scale=0.33]{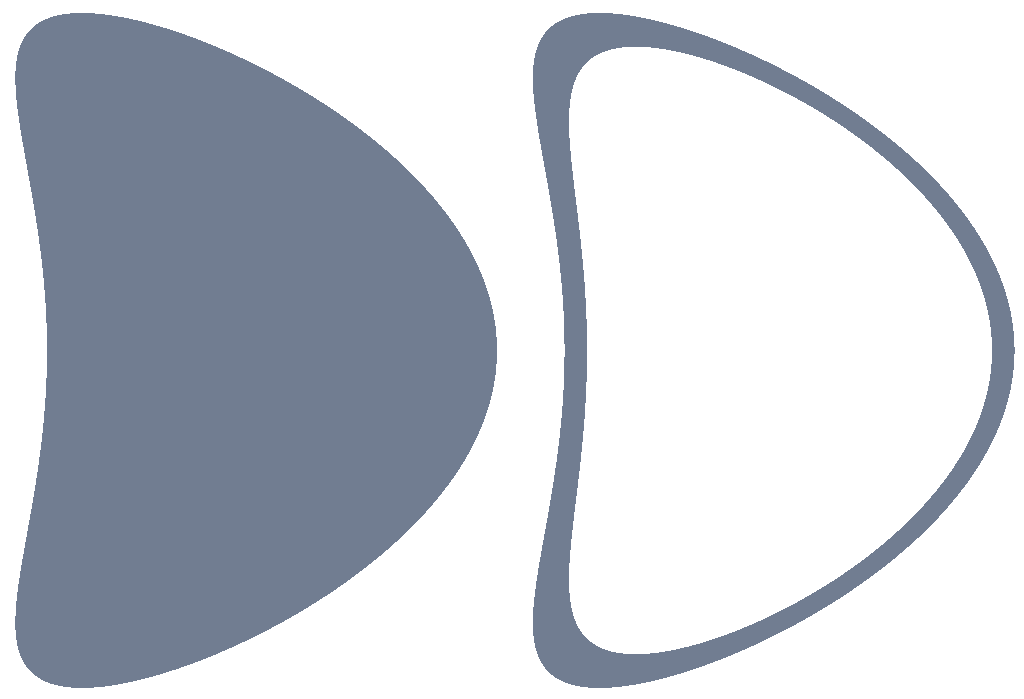}}	
		\hfill	
		\subfigure[A square base region $\Omega_{E}$ containing $\Omega$.]	{\includegraphics[clip=true, trim=000 00  00 00, scale=0.3]{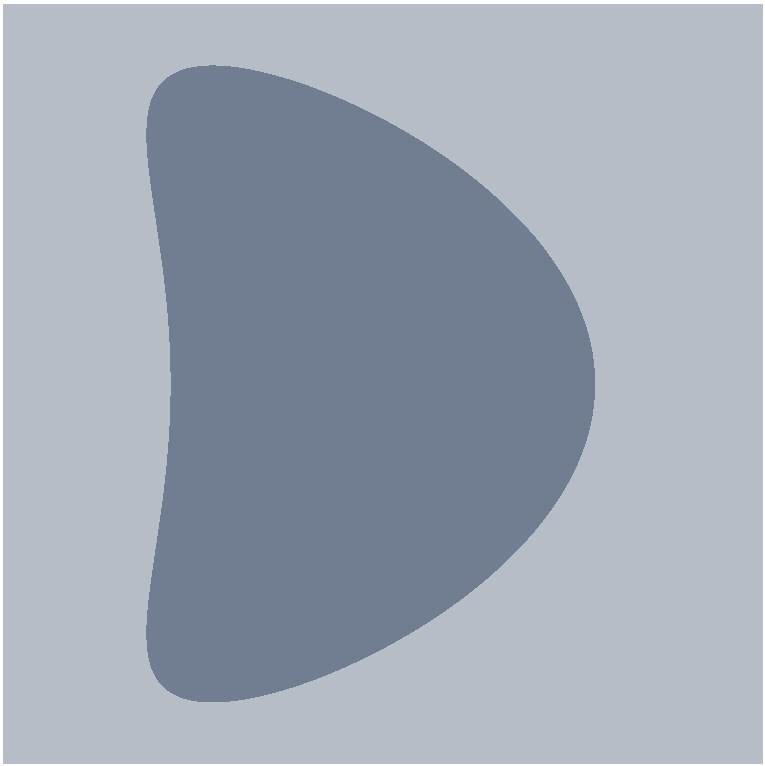}}		
		\hfill
		\subfigure[A cut-off function $\eta_{\tau_0}((\mcP\bsx-\bsx)\cdot \bsnu(\mcP\bsx))$ in the base region $\Omega_{E}$.]	{\includegraphics[clip=true, trim=000 00  00 00, scale=0.40]{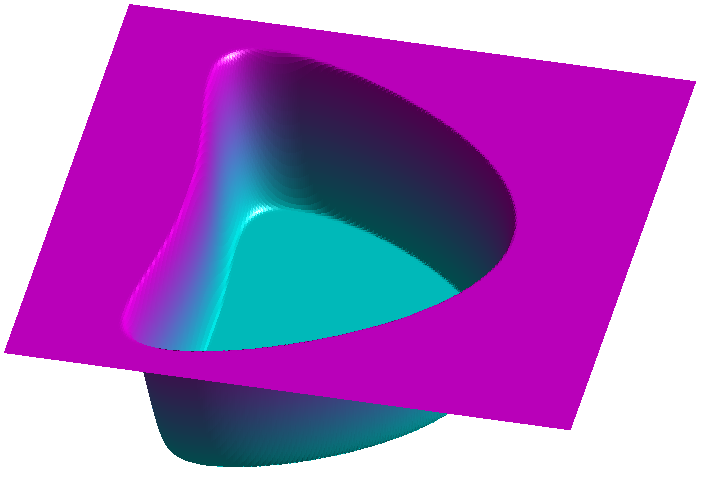}}	
		\hfill	
		\subfigure[The corresponding $1-\eta_{\tau_0}((\mcP\bsx-\bsx)\cdot \bsnu(\mcP\bsx))$ on base region $\Omega_E$.]	{\includegraphics[clip=true, trim=00 0  140 00, scale=0.17]{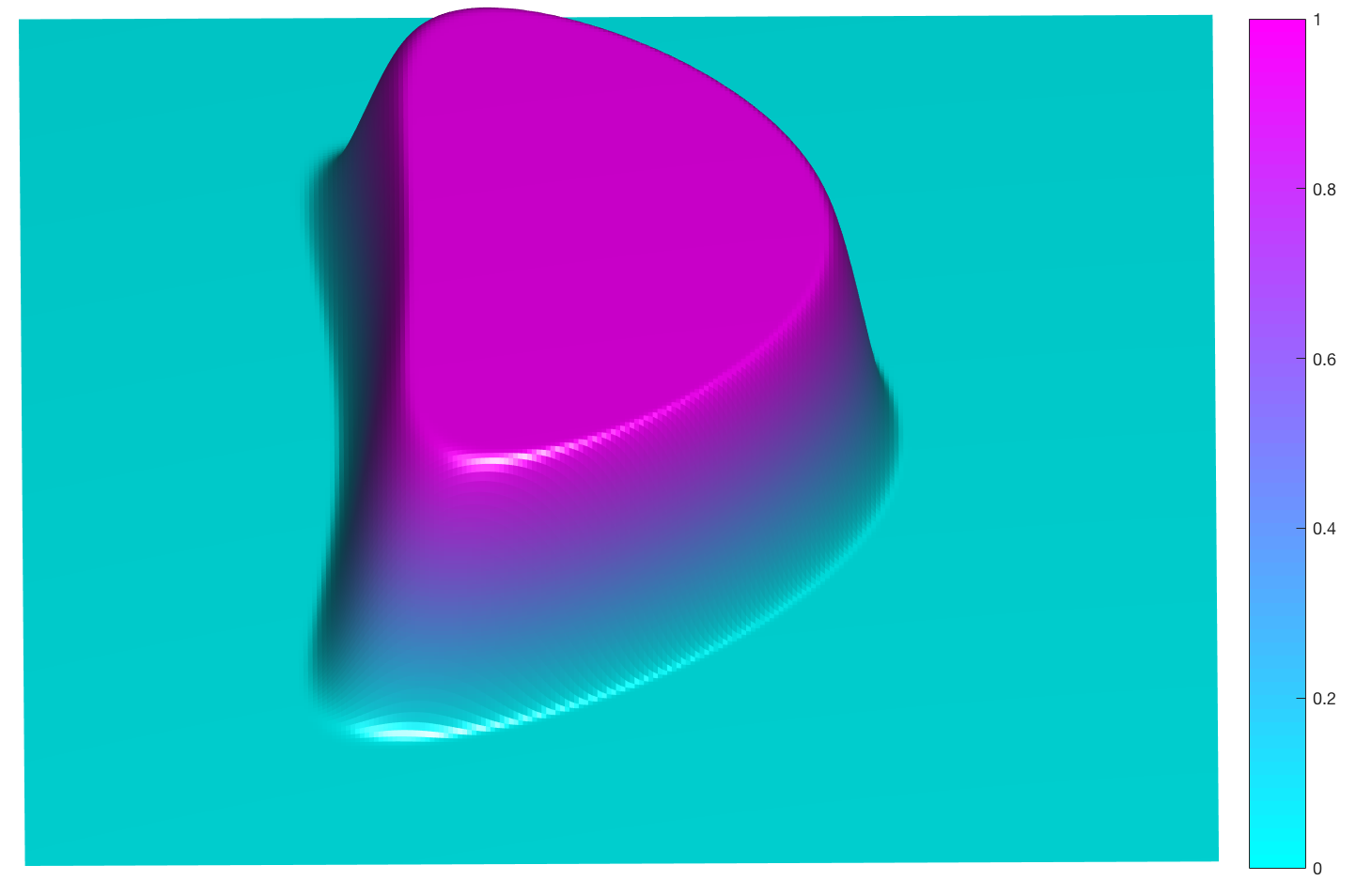}}	
		
		\caption{A splitting of $\Omega$ into boundary ($\Omega_{B}$) and base ($\Omega_{E}$) regions  with corresponding cut-off functions. In the region encompassed by $\Omega_B$, $1-\eta_{\tau_0}$ takes value one, transitions smoothly to zero as we approach the boundary of $\Omega$ and stays there in $\Omega_{E} \setminus \Omega$.}		\label{fig:-omega_spliting}		
	\end{center}
\end{figure}

\subsection{Boundary and base integrals}
To isolate the boundary region for a specialized treatment, we utilize the cut-off function
\begin{equation}\label{eq:-cut-off}
\eta_{\tau_0} (\tau) = 
\begin{cases}
1, \hspace{1mm} &\text{for }  \tau \le 0,\\
 \exp\left( \frac{2 \tau_0 e^{-\tau_0/\tau}}{\tau-\tau_0} \right),&  \text{for } 0<\tau <\tau_{0}, \hspace{1mm}\\
 0,& \text{for } \tau_{0} \le \tau,
\end{cases}
\end{equation}
(shown in Figure \ref{fig:-cutoff}) to decompose $\mcA(v)$ as
\begin{equation} \label{eq:-boundary-bulk}
\mcA(v)(\bsx) = \int \limits_{\Omega_B}G_{\kappa}(\bsx - \bsy) v(\bsy) \eta_{\tau_0}(\tau')d\bsy  + \int \limits_{\Omega}G_{\kappa}(\bsx-\bsy)v(\bsy) (1-\eta_{\tau_0}(\tau'))d\bsy 
\end{equation}
where $\tau' = (\mcP\bsy-\bsy)\cdot \bsnu(\mcP\bsy)$. Note that, owing to presence  of the factor $1-\eta_{\tau_{0}}(\tau')$, the domain of integration of the second integral in (\ref{eq:-boundary-bulk}) can be extended arbitrarily, to say $\Omega_E$, by extending the integrand by zero without effecting its smoothness. In particular, in this framework, we always choose $\Omega_{E}$ in such away that $\Omega\subset \Omega_{E}$. A specific instance of this procedure, for the case of bean shape scatterer, is displayed in Figure \ref{fig:-omega_spliting}. To make this more precise, we define $\mcA_E$ by
\begin{equation}
\label{eq:-AInt}
\mcA_E (v)(\bsx) = \int_{\Omega_E} G_{\kappa}(\bsx - \bsy) \mcE(v)(\bsy) d\bsy 
\end{equation}
where
\begin{equation}
\mcE(v)(\bsx) = 
\begin{cases}
v(\bsx), & \bsx \in \Omega \setminus \Omega_B, \\
v(\bsx) \left(1-\eta_{\tau_0}((\mcP\bsx-\bsx)\cdot \bsnu(\mcP\bsx))\right), & \bsx \in \Omega_B, \\
0, & \bsx \in \Omega_E \setminus \Omega, \\
\end{cases}
\end{equation}
and denote the first integral in (\ref{eq:-boundary-bulk}) as
\begin{equation}
\label{eq:-ABdry}
\mcA_B (v)(\bsx) = \int_{\Omega_B} G_{\kappa}(\bsx - \bsy) v(\bsy) \eta_{\tau_0}((\mcP\bsy-\bsy)\cdot \bsnu(\mcP\bsy)) d\bsy,
\end{equation}
thus arriving at the final boundary-base decomposition of $\mcA(v)$ given by
\begin{equation}
\label{eq:-split}
\mcA(v)(\bsx) = \mcA_B(v)(\bsx) + \mcA_E(v)(\bsx).
\end{equation}

\subsection{Discretization and grid interactions}
\label{sec:-Dis_int}
At this stage, we present an overview  of our discretized problem. We make these details more precise in subsequent sections where we apply our framework to produce specific instances of scattering solvers. To begin with, our overall computational grid, say $\Omega_h^o$, comprise of meshes on boundary and base regions that we denote by  $\Omega_{B,h}$ and $\Omega_{E,h}$ respectively. While the nodes used for $\Omega_{E,h}$ depend on the fast convolution technique being used for base integration, the grid on $\Omega_B$ typically follows a parametric description of the boundary region that involves a set of overlapping subregions, each homeomorphic to $[0,1]^2$. In fact, the nodes in $ \Omega_{B,h}$ in each of these sub-boundary regions are the direct image of a regular mesh on $[0,1]^2$ under the corresponding parametric map.

\begin{figure}[t] 
	\begin{center}	
		\subfigure	{\includegraphics[clip=true, trim=13 10 20 8, scale=0.55]{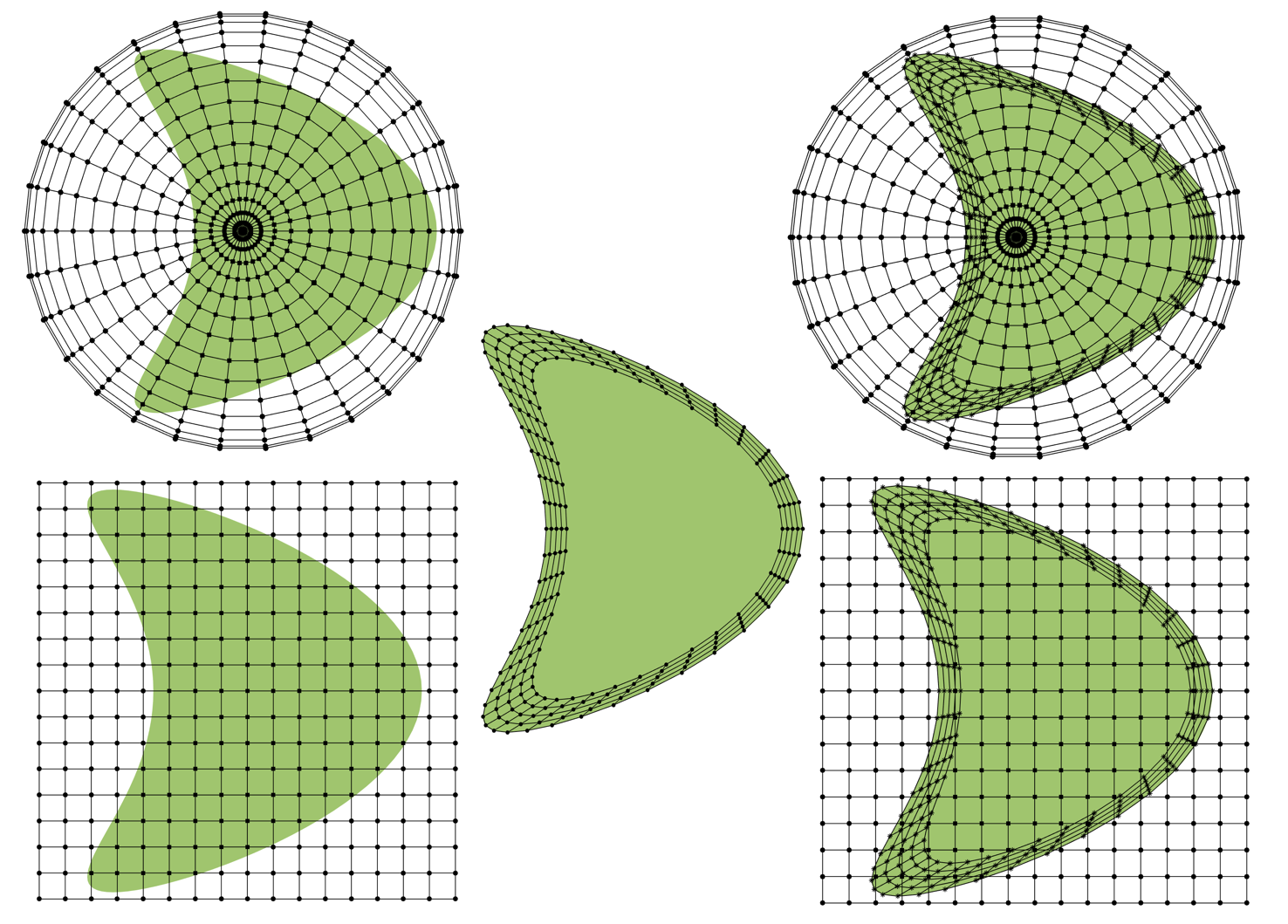}}
		\caption{
			Example computational grids $\Omega_{h}$ for proposed schemes with the Addition theorem approach and the pre-corrected trapezoidal method. 
			A polar base grid $\Omega_{E,h}$ required for approximation of base integral using the Addition theorem method is shown on {\bf top right} whereas the figure in the {\bf center} depicts a boundary grid $\Omega_{B,h}$. While the union of base and boundary grids form the overall computational grid $\Omega_h^o$, the solution is computed only on the subgrid $\Omega_h$ contained in $\Omega$. The corresponding combined view of boundary and base computational grids is shown on {\bf top right}. 
			The figures on {\bf bottom left} and {\bf bottom right} similarly show grids for the pre-corrected trapezoidal method under our framework. 
		}
					\label{fig:-grid}		
	\end{center}
\end{figure}

While we do use $\Omega_h^o$ for the computation of $\mathcal{A}(v)$, as necessitated by requirements of underlying fast convolution techniques, only a subset, that we refer to as approximation grid and denote by $\Omega_h$, is used as nodes for discretization of the integral equation \ref{eq:-Lippmann}. To be more precise, if set of grid points in $\Omega_{E,h}$ which do not lie in $\Omega$ is denoted by $\Omega_{E,h}^{e}$, then
\[
\Omega_h = \Omega_{B,h} \cup \left(\Omega_{E,h} \setminus \Omega_{E,h}^{e}\right).
\]
 A couple of example computational grids are shown in Figure~\ref{fig:-grid} for illustration.

In view of (\ref{eq:-split}), the evaluation of approximations $\mcA_h(v) = \mcA_{B,h}(v) + \mcA_{E,h}(v)$ for points in $\Omega_h$, can be split into four distinct calculations, namely, 
\begin{enumerate}
\item $\mcA_{E,h}(v)(\bsx), \bsx \in  \Omega_{E,h}$, the base-base interaction,
\item $\mcA_{B,h}(v)(\bsx), \bsx \in  \Omega_{B,h}$, the boundary-boundary interaction,
\item $\mcA_{E,h}(v)(\bsx), \bsx \in  \Omega_{B,h}$, the boundary-base interaction, and 
\item $\mcA_{B,h}(v)(\bsx), \bsx \in  \Omega_{E,h} \setminus \Omega_{E,h}^{e}$, the base-boundary interaction.
\end{enumerate}

%

We provide details for each of these interactions in the following four sections.

\section{Base-base interaction}
\label{sec:-BubI}
As $\Omega_E$ covers the base of the scattering media, calculations pertaining to the base-base interactions are most voluminous among all interactions and, therefore, demands a special attention with respect to the speed of computation, of course, without undue compromise in the accuracy of approximations. 
To this end, we recall that the factor $\mcE(v)$ in the integrand of  $\mcA_E(v)$ vanishes to high order on the boundary of $\Omega_E$, and therefore, any fast high-order numerical integration scheme that converges rapidly while computing the convolution of a singular kernel and a globally smooth function can be employed for accurate calculation of this interaction. To demonstrate that this indeed can be achieved, and to supply underlying details, we have used two entirely different approaches, namely, a pre-corrected trapezoidal rule \cite{aguilar2004high,duan2009high}  and an Addition theorem based method \cite{bruno2004efficient} for performing these calculations. We succinctly review both this approaches in the following subsections.

\subsection{Base integration with pre-corrected trapezoidal rule}
It is well known that  trapezoidal rule for smooth, periodic integrands converges super-algebraically fast. However, it does not yield satisfactory results when used to compute convolutions involving singular kernels. To overcome this difficulty, one can utilize a pre-corrected trapezoidal rule, originally proposed by Rokhlin in \cite{rokhlin1990end} for a certain class of singular functions which has subsequently been developed further by several researchers. 
 The pre-corrected trapezoidal rule consists of classical trapezoidal rule punctured at the point of integrand singularity along with a correction operator incorporated in the vicinity of singularity. A weighted sum of integrand values in a neighborhood of singular point constitutes the correction operator where the newly introduced weights are chosen in a manner that alleviates the singular behavior of the integrand to achieve high-order convergence. Further, the correction term does not require evaluation of singular function at the point of singularity. Additionally, these weights can be pre-computed by solving a linear system and can be tabulated for further use. We must note, however, that large condition numbers of matrices being inverted during this process require that these linear systems be solved to a very high precision using multi-precision computer arithmetic. 
The two dimensional pre-corrected trapezoidal  quadrature for the case of logarithmic singularity  appeared in \cite{duan2009high,marin2014corrected,aguilar2002high}.   The method presented in \cite{duan2009high},  provides an analytical expression for the calculation of correction weights and thereby avoids solving ill-conditioned linear systems. While the analytical calculation of these formulas are non-trivial, they can be evaluated numerically to a high accuracy for further use. However, as these numerical computations are tied to the underlying mesh, these computationally large pre-computations, in practice, may be required for every scattering simulation.
On the other hand, correction weights computed in  \cite{marin2014corrected,aguilar2002high}, are obtained as solutions to a linear system and present an attractive alternative to the one mentioned above. In fact, in our implementations, to avoid the numerical computation of correction weights, we have used those that are tabulated in \cite{aguilar2002high}. 

 As mentioned in the beginning of this section, $\mcE(v)$ is smooth and compactly supported within $\Omega$, and therefore, the domain of integration in (\ref{eq:-AInt}) can be extended to an arbitrary $\Omega_E \supset \Omega$ without altering its value. 
For employing the pre-corrected trapezoidal rule, we use $\Omega_E = [-a,a]\times [-a,a]$,  where the positive real number $a$ is chosen large enough so that $\Omega \subset \Omega_E$.  We then set a computational grid $\Omega_{E,h}$ on $\Omega_{E}$ given by 
 \begin{equation} \label{eq:-omegaeh}
 	\Omega_{E,h} =\{(ih,jh)| \hspace{1mm} -n \leq i,j \leq n, n=a/h  \}, 
 \end{equation}
 where $h$ denotes the mesh size.  
 Next, we partition  $\Omega_{E,h}$ into disjoint sets according their 
distance from the origin. 
To be precise, any two grid points   $\bsx_{1}=\left(i_1h,j_1h \right),
\bsx_{2}=\left(i_2h,j_2h \right) \in \Omega_{E,h}$
belong to the same set if $i_1^2+j_1^2 =
i_2^2+j_2^2$. In each set, say $S_{r}$, there exist unique grid point $\tilde{\bsx}_{r}=(ih,jh)$ such that 
$0\leq j \leq i$, where the subscript $r$ is given by $r=i(i+1)/2+j+1$. 

 We re-expressed integral kernel $G_{\kappa}(\bsx - \bsy) $  as a sum of non-smooth and smooth function as \cite{kress1989linear},
\begin{equation} \label{eq:-kernelsplit}
 G_{\kappa}\left(\bsx - \bsy\right) =P_{\kappa}\left(\bsx - \bsy \right)\ln \left(|\bsx-\bsy|\right)+Q_{\kappa}\left(\bsx - \bsy\right), 
\end{equation}
 where 
 \begin{align*}
 P_{\kappa}\left(\bsx\right) &= -\frac{J_0\left(\kappa|\bsx |\right)}{2\pi},\\
 Q_{\kappa}\left(\bsx\right) &=  -\frac{1}{2\pi} \left[ \ln\left(\frac{\kappa}{2} \right)+\gamma\right]+\frac{i}{4}+R_{\kappa}\left(|\bsx|\right),\\
 \gamma &= 0.5772156649015328606 \cdots
 \end{align*}
 is Euler's constant, and $R_{\kappa}$ being smooth function such that $R_{\kappa}(0) =0$.
  Note that, for any point 
 $\bsx\in S_r$, $P_{\kappa}\left(\bsx\right)=P_{\kappa}\left(\tilde{\bsx}_{r}\right)$. 
 Now,  $\mcA_E (v)$ can be expressed as
 \begin{equation}\label{eq:-BBI-1}
 \mcA_E (v)(\bsx) = \int_{\Omega_E} P_{\kappa}\left(\bsx - \bsy \right) \ln \left(|\bsx-\bsy | \right)\mcE(v)(\bsy) d\bsy +\int_{\Omega_E} Q_{\kappa}\left(\bsx - \bsy\right) \mcE(v)(\bsy) d\bsy.
 \end{equation}
  In the above equation, integrand of the second integral on the right hand side  is smooth and vanishes toward the boundary of integration region, therefore, high-order approximations can be achieved by means of classical trapezoidal rule. On the other hand, integrand in the first integral has a logarithmic singularity at $\bsx =\bsy$, which we integrate to high-order by using the corresponding pre-corrected trapezoidal quadrature.
  If $\bm{c}_{r}, r=1,\cdots,k$, denotes the $k$ correction coefficients, then for any grid point $\bsx_{l}\in \Omega_{E,h}$, application of pre-corrected trapezoidal rule on first integral and classical trapezoidal rule on second integral of (\ref{eq:-BBI-1}) gives 
  \begin{align}\label{eq:-BBI-2}
  \mcA_E (v)(\bsx_{l}) &= h^{2}\sum_{\substack{
  		\bsx_{j}\in \Omega_{E,h} \\ \nonumber
  		\bsx_{l} \neq\bsx_{j}
  	}} 
  	P_{\kappa}\left(\bsx_{l}-\bsx_{j}\right) \ln\left(|\bsx_{l}-\bsx_{j} |\right)\mcE(v)(\bsx_{j})
  	+h^2\sum_{	\bsx_{j}\in \Omega_{E,h}}	Q_{\kappa}\left(\bsx_{l}-\bsx_{j}\right) \mcE(v)(\bsx_{j})  	
  \\ \nonumber
  	&+	h^2 \ln(h)
  	P_{\kappa}\left(\bm{0}\right)\mcE(v)(\bsx_{l})
  		+h^2\sum_{r=1}^{k}c_{r} d_{r}\sum_{\bsx_{q}\in S_{r}}\mcE(v)\left(\bsx_l+\bsx_{q}\right)+O\left(h^{2p+4} \right),
 \end{align}
 where $d_r= P_{\kappa}\left(\tilde{\bsx}_{r}\right)$ and $p$ is positive integer related to the number of correction coefficients $k$ according to $k=p(p+1)/2+1$. As, 
  \begin{equation}
  Q_{\kappa}\left(\bm{0}\right) =  -\frac{1}{2\pi} \left[ \ln\left(\frac{\kappa}{2} \right)+\gamma\right]+\frac{i}{4},
  \end{equation}
 therefore, above equation can be simplified as 
  \begin{equation}\label{eq:-BBI-3}
  \mcA_E (v)(\bsx_l) = h^{2}\sum_{\substack{
  		\bsx_{j}\in \Omega_{E,h} \\ 
  		\bsx_l \neq\bsx_{j}
  	}} 
  	G_{\kappa}(\bsx_l - \bsx_j)\mcE(v)(\bsx_{j})	 
  	+h^2\sum_{r=1}^{k}w_r\sum_{\bsx_{q}\in S_{r}}\mcE(v)(\bsx_l+\bsx_{q})+O\left(h^{2p+4} \right),
  	\end{equation}
  	where
  	\begin{align*}  	
  	w_1 &= -\frac{1}{2\pi} \left(\ln\left( h\kappa/2\right) +c_1+\gamma-\frac{2\pi i}{4}\right),\\ 
  	w_r &=d_rc_r, r=2,\cdots,k. 
  	\end{align*}  
  Note that, in Eq. (\ref{eq:-BBI-3}), terms in the summation over $x_q \in S_r$ are zero whenever $x_l+x_q \not\in \Omega_{E,h}$.
  Thus, for any positive  $ p$, we can approximate integral operator  $\mcA_{E,h}(v)(\bsx) $ at all grid points with accuracy of order $O\left(h^{2p+4}\right)$ by correcting the integral weight at $k=p(p+1)/2+1$ grid points near the vicinity of singular point. Further, as grid points are equidistant, the discrete convolution in  (\ref{eq:-BBI-3}) can be obtained in $O(N_{I}\log N_I)$ operations by means of FFT, where $N_{I}=4n^2$ denotes the total number of base grid points.     

\subsection{Base integration with the Addition theorem approach}

In this case, we choose $\Omega_E$ to be a disc of radius $R$ so that $\Omega \subset \Omega_E$. In polar coordinates with $\bsx = (r\cos\theta,r\sin\theta)$ and $\bsy = (r'\cos\theta',r'\sin\theta')$, using the Addition theorem for the Hankel function \cite{colton2012inverse},  (\ref{eq:-AInt}) is rewritten as
\begin{equation}
\label{eq:-AIntPolar}
\mcA_E(v)(\bsx) = \sum_{\ell = -\infty}^{\infty} e^{i\ell \theta} \int_{0}^R G_{\kappa,\ell}(r,r') \left[  \int_0^{2\pi} \mcE(v)(\bsy)e^{-i\ell\theta'}\,d\theta' \right] \,dr',
\end{equation}
where, denoting the first kind Bessel and Hankel functions of order $\ell$ by $J_{\ell}$ and $H_{\ell}^1$ respectively,
\[
G_{\kappa,\ell}(r,r') =  \frac{i}{4} H_{\ell}^1(\kappa\max(r,r'))J_{\ell}(\kappa\min(r,r')).
\]
For any fixed $r \in [0,R]$, the evaluation of $\theta'$-integral in (\ref{eq:-AIntPolar}) amounts to computing the $\ell^{th}$ Fourier coefficient of a $2\pi$ periodic function, and can be affected to high-order accuracy, for example, through the use of trapezoidal rule. The integrand in the radial integration, however, is only piecewise smooth in intervals $(0,r)$ and $(r,R)$. 
Denoting the $\ell^{th}$ Fourier coefficient of $g$ by $(g)_{\ell}$,  the second kind Bessel function of order $\ell$ by $Y_{\ell}$, and following \cite{bruno2004efficient}, a high order accuracy in the approximation is obtained by rewriting (\ref{eq:-AIntPolar}) as
\[
\mcA_E(v)(\bsx) = -\frac{i\pi}{2} \sum_{\ell = -\infty}^{\infty} e^{i\ell \theta} \left[  F^{(1)}(r) + F^{(2)}(r) -i \frac{J_{\ell}(\kappa r)}{Y_{\ell}(\kappa R)}F^{(1)}(R) \right]
\]
where
\begin{equation}
F^{(1)}(r) = \int_0^r Y_{\ell}(\kappa r) J_{\ell}(\kappa r') (\mcE(v))_{\ell}(r') \,dr'
\end{equation}
and
\begin{equation}
F^{(2)}(r) = \int_r^R J_{\ell}(\kappa r) Y_{\ell}(\kappa r') (\mcE(v))_{\ell}(r') \,dr'.
\end{equation}
Though we've uniform angular grid for computing $(\mcE(v))_{\ell}(r')$, our radial grid for this interaction is somewhat different and will be explained when we discuss evaluation procedure of $F^{(1)}(r)$ and $F^{(2)}(r)$. Discontinuities and corner singularities in $F^{(1)}(r)$ and $F^{(2)}(r)$ due to $(\mcE(v))_{\ell}(r')$ are resolved by breaking the (radial-) integration domain $[0,R]$ into several subintervals of uniform length, where the breaking points are inclusive of the aforementioned singular points. Further, in each of those subintervals, we have Chebyshev grid of order $N_c$ and for rest of this section we concentrate on one of these uniform singularity free subintervals. Let $[a,b]$ be a such interval and $\alpha_{ab}(r')=\frac{r'-a}{b-a}-\frac{b-r'}{b-a}$ maps $[a,b]$ onto the standard interval $[-1,1]$. we do an approximation of 
$(\mcE(v))_{\ell}(r')$ via a truncated Chebyshev series 
$$(\mcE(v))_{\ell}(r')=\displaystyle\sum_{n=0}^{N_c-1}c_nT_n(\alpha_{ab}(r')),$$
which provides approximation almost as good as minimax polynomial. 
For $\beta\in [-1,1]$, the inverse of $\alpha_{ab}$ is given by
$$\alpha_{ab}^{-1}(\beta) = \frac{1}{2}\big[(b-a)\beta+b+a \big].$$
We thus need the values of $(\mcE(v))_{\ell}(r')$ at the points 
$a_j = \alpha_{ab}^{-1}(\beta_j)$, where $\beta_j = \cos\big(\pi(j-0.5)/N_c\big)$ for $j=1,2,..,N_c$ are corresponding Chebyshev points in the standard interval $[-1,1]$, to compute the coefficients $c_n$. We also include the endpoints in the radial grid and denote them, for convenience, as $a_0 = a$ and $a_{N_c+1 } = b$. 
We then break the integrals over $[a,b]$ as follows
\begin{equation}
\int_{a_0}^{a_{N_c+1}} Y_{\ell}(\kappa r) J_{\ell}(\kappa r') (\mcE(v))_{\ell}(r') \,dr' = \displaystyle\sum_{j=0}^{N_c}c_j\int_{a_j}^{a_{j+1}}Y_\ell(\kappa a_j)J_\ell(\kappa r)T_n(\alpha_{ab}(r'))r'dr' 
\end{equation}
\begin{equation}
\int_{a_0}^{a_{N_c+1}} J_{\ell}(\kappa r) Y_{\ell}(\kappa r') (\mcE(v))_{\ell}(r') \,dr' = \displaystyle\sum_{j=0}^{N_c}c_j\int_{a_j}^{a_{j+1}}J_\ell(\kappa a_j)Y_\ell(\kappa r)T_n(\alpha_{ab}(r'))r'dr'
\end{equation}
that allows for one time computation of the integral moments
\begin{equation}
  P^{(1)}_{\ell n} = \displaystyle\int_{a_0}^{a_{j}}J_\ell(\kappa a_j)Y_\ell(\kappa r)T_n(\alpha_{ab}(r'))r'dr'
\end{equation}
\begin{equation}
  P^{(2)}_{\ell n} = \displaystyle\int_{a_j}^{a_{N_c+1}}J_\ell(\kappa a_j)Y_\ell(\kappa r)T_n(\alpha_{ab}(r'))r'dr'
\end{equation} 
that can be stored for repeated use in further calculations. An accurate approximation of the radial integral, of course, can be obtained from these moments by adding and subtracting scaled values of $P^{(1)}_{\ell n}$ and $P^{(2)}_{\ell n}$, as necessary, in $\mathcal{O}(N)$ operations.

\section{Boundary-boundary interaction}
\label{sec:-BBI}

\begin{figure}[h!] 
\begin{center}	
{\includegraphics[clip=true, trim=0 0  0 0, scale=0.67]{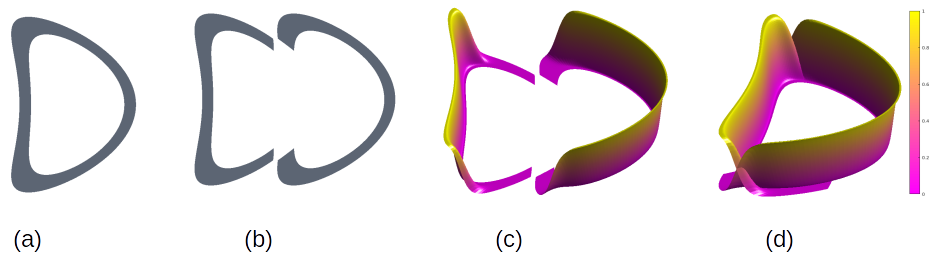}}	
\caption{For a bean shape geometry, boundary region, overlapping coordinate patches and windowing function. Boundary region in (a) represented as two overlapping coordinate patches are shown in (b), corresponding windowing functions $w_{k}\eta_{\tau_{0}}$  are depicted in (c) while  (d) displays a combined view.}
\label{fig:-boundarypou}
\end{center}
\end{figure}

We begin with a $K$-piece overlapping cover of the boundary curve $\partial \Omega$, say $\{ \partial\Omega_k \}_{k = 1}^K$, each of which has a smooth invertible parametrization  $\psi_k : [0,1]\to \partial\Omega_k$. A parametrization $\bsxi_k = \bsxi_k(s,t) : [0,1]^2 \to \Omega_{k,B}$ for the boundary sub-region $\Omega_{k,B}$, defined by
\[
\Omega_{k,B} = \{ \bsx \in \Omega_B \ :\  \mcP\bsx \in \partial\Omega_k \},
\]
is obtained from $\psi_k$ where its inverse is given by
\[
\bsxi_k^{-1}(\bsx) = (s,t)=(\psi^{-1}_{k}(\mcP\bsx),(\mcP\bsx-\bsx)\cdot \bsnu(\mcP\bsx)/\tau_0).
\]
The smoothness of $\bsxi_{k}$ follows from the smoothness of $\partial\Omega$.
Now, with the help of a
partitions of unity $\left\{w_k :k = 1,...,K\right\}$ subordinate to the covering $\{ \Omega_{k,B} \}_{k = 1}^K$ of $\Omega_B$, specific instances of which have been shown in Figure \ref{fig:-boundarypou}, the integral in (\ref{eq:-ABdry}) is rewritten as a sum of integrals,
\begin{equation} 
 \label{eq:-BrdySplit}
  \mcA_B(v)(\bsx)= \sum_{k=1}^{K} \mcA_{k,B}(v)(\bsx)
\end{equation}
  where
 \begin{equation} 
 \label{eq:-BdrySplitK}
  \mcA_{k,B}(v)(\bsx) =\int_{0}^{1} \int_{0}^{1} G_{\kappa}(\bsx-\bsxi_k(s',t')) v(\bsxi_{k}(s',t')) \eta_{\tau_0}(\tau_0 t')\xi^{\prime}_k(s',t') w_k(\bsxi_k(s',t'))\,ds' \,dt'.
 \end{equation} 
Here, $\xi^{\prime}_k$ is the Jacobian of the transformation $\bsxi_k$. 
Although,  the boundary region $\Omega_B$ is typically a small fraction of the scattering medium $\Omega$, nonetheless, a brute force integration requires  $O\left(NN_{B}\right)$ operations to evaluate $ \mcA_B(v)(\bsx)$ for all $\bsx \in \Omega_{h}$, where $N_{B}$ denote the cardinality of the set $\Omega_{B,h}$.

In order to reduce this computational cost for the overall discretization scheme, we split  computation in (\ref{eq:-BrdySplit}) into two, namely, adjacent and non-adjacent interactions. To give the precise definition of adjacency, we introduce a square cell $\mathcal{C}$ of side length $A$, containing the scatterer $\Omega$.  We then partition the square $\mathcal {C}$ into $L^2$ identical cell $c_{ij}$ $\left(i,j=1,\cdots L \right)$ of side length $H=A/L$, such that there are $L$ cell along each side of the square.  For each discretization point $\bsx \in c_{i_{0}j_{0}}$ which we often refer as a source point, we define adjacent set $\mathcal{N}_{\bsx}$ as

\begin{equation}
\label{eq:-adj}
\mathcal{N}_{\bsx}=\left\{ \bsy \in c_{ij} \ | \hspace{2mm} | i-i_{0}| \leq 1, | j-j_{0}| \leq 1\right\}.
\end{equation}

Note that, for each $\bsx \in c_{i_0j_0}$, $\mathcal{N}_{\bsx}$ contains at most eight neighboring cells.
A source point $\bsy \in c_{ij}$ is said to be adjacent to $\bsx$ if $\bsy \in \mathcal{N}_{\bsx}$ and non-adjacent if otherwise.  A specific illustration of this procedure displayed in Figure \ref{fig:-accelpic}. Using the cut-off function $\eta_{s}(s')$, supported within adjacent set $\mathcal{N}_{\bsx}$ and choosing the thickness of boundary region smaller than the $3H$, evaluation of integral (\ref{eq:-BdrySplitK}) can by decomposed in adjacent and non-adjacent interactions as follows: 
 \begin{equation} 
\label{eq:-BdrySplitK2}
\mcA_{k,B}(v)(\bsx) =\int_{0}^{1} \int_{0}^{1} \cdots \eta_{s}(s')w_k(\bsxi_k(s',t'))\,ds' \,dt'+
\int_{0}^{1} \int_{0}^{1} \cdots (1-\eta_{s}(s'))w_k(\bsxi_k(s',t')) \,ds' \,dt'.
\end{equation}
Note that, the second integral is always non-singular while the first is singular if target point $\bsx\in \Omega_{k,B}$. In order to achieve desirable accuracy while reducing overall cost of the algorithm, we employ different strategies for evaluation of adjacent (singular) and non-adjacent (nonsingular) integrals, which we discuss next.

\subsection{Adjacent interaction }
\label{subsec:-adjI}
 When $\bsx \not\in \Omega_{k,B}$, adjacent integral  can be handled by use of high-order quadrature rules in both $s'$ and $t'$ variables as integrand is smooth. 
We observe that the 
$w_k(\bsxi_k(s,t))$ vanish to high-order at $s' = 0$ and $s' = 1$. On the other hand, while $\eta_{\tau_0}(\tau_{0}t')$ vanishes at $t' = 0$, it does not vanish at $t' = 1$ (for instance see Figure \ref{fig:-boundarypou}). Thus, for $s'$-integration, trapezoidal rule converges rapidly due to high-order vanishing integrands, while in $t'$ direction, a high-order composite Newton-Cotes quadrature is employed to produce accurate approximations. This, of course, requires that a uniform grid be placed in the $s-t$ parameter space, $[0,1]^2$. In our implementation, we have used five point composite Newton-Cotes quadrature for the transverse $t'$-integration. 

If $\bsx \in \Omega_{k,B}$, say $\bsx = \bsxi_k(s,t)$, then the integrand  of first integral in (\ref{eq:-BdrySplitK2}) is clearly singular when $(s',t') = (s,t)$. Moreover, even when $t' \ne t$, the kernel  $G_{\kappa}(\bsxi_k(s,t)-\bsxi_k(s',t'))$ as a function of $s'$, though not singular, exhibit increasingly rapid variation as $t'$ approaches $t$. Following \cite{anand2006efficient}, a change of variable $s' = s'(\tau)$ with $s'(\tau) = s + \rho(\tau)$, $\rho$  being a smooth invertible odd function such that
\begin{equation}
\left.\frac{d^m\rho}{d\tau^m}\right|_{\tau = 0} = 0, \text{\ \ \ for } m = 0,\ldots,M,
\end{equation}
is introduced to overcome both difficulties. Indeed, as shown in \cite{anand2006efficient}, the integral
\begin{equation}
\label{eq::transverse}
I^k_{s,t}(t') = \int\limits_{-\rho^{-1}(s)}^{\rho^{-1}(1-s)} G_{\kappa}(\bsxi_k(s,t)-\bsxi_k(s'(\tau),t')) \varphi(\tau,t') \rho^{\prime}(\tau) \,d\tau  
\end{equation}
with
$$\varphi(\tau,t') = v(\bsxi_k(s'(\tau),t')) \eta_{s}(s')w_k(\bsxi_k(s'(\tau),t')) \xi_k(s'(\tau),t')$$ 
can be evaluated to high-order using a trapezoidal rule in the variable $\tau$. Evaluation of $\varphi$ on a uniform $\tau$-grid, however, requires the $v$ data at off-grid points that necessitates use of an efficient and high-order accurate interpolation scheme. In this connection, a mixed trigonometric-polynomial interpolator, introduced in \cite{bruno2001fast}, is used, which is known to be an effective strategy when the underlying function has a smooth periodic extension.
As shown in \cite{anand2006efficient}, the $t'$-integrand is piecewise smooth in $(0,t)$ and $(t,1)$ and therefore can be integrated to high-order by first breaking the transverse integral in (\ref{eq:-BdrySplitK2}) at $t' = t$ as
\[
\int_{0}^{1} \int_{0}^{1} \cdots \eta_{s}(s')w_k(\bsxi_k(s',t'))\,ds' \,dt' = \int_0^1 I^k_{s,t}(t')\eta_{\tau_0}(\tau_0 t')\,dt' =  \int_0^t I^k_{s,t}(t')\eta_{\tau_0}(\tau_0 t')\,dt' +  \int_t^1 I^k_{s,t}(t')\eta_{\tau_0}(\tau_0 t')\,dt'
\]
and then using a high-order composite Newton-Cotes quadrature for the approximation of each of the integrals. However, as explained in \cite{Anand2015highorder}, this strategy requires that we know $I^k_{s,t}$ at off-grid points, particularly near $t' = 0$ and $t'=1$. Clearly, the direct interpolation of $I^k_{s,t}$ is not high-order accurate in view of the corner singularity at $t'=t$. Again, following the strategy in \cite{Anand2015highorder}, we resolve this by adding $2\times(Q-1)\times(Q-1)$ additional grid, $(Q-1)\times(Q-1)$ in the vicinity of zero and $(Q-1)\times(Q-1)$ in the vicinity of one, where $Q$ is the order of Newton-Cotes quadrature rule.  The values of $I^k_{s,t}$ at these additional grid points are obtained by, first interpolating smooth density $\varphi(\tau,t')$ at these extra grid points followed by integration in  (\ref{eq::transverse}).
In light of the fact that, each nonempty cell $c_{ij}$ contains on an average $O(N_{B}/L^2)$ discretization points, this procedure for evaluation of adjacent interactions require $O(N_B^2/L^2)$ operations in total.

\subsection{Non-adjacent interaction}
\label{sec:-nsnai}
\begin{figure}[h!]
	\begin{center}			
		{\includegraphics[clip=true, trim=0 0  00 0, scale=0.7]{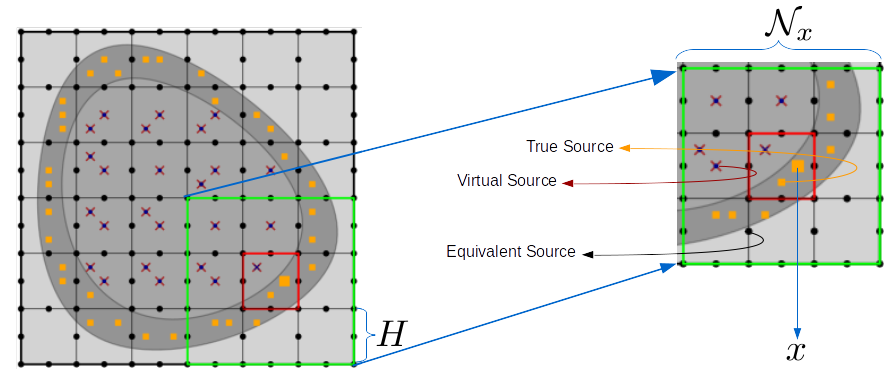}}	
		\caption{An illustration of the non-singular non-adjacent boundary interaction by means of two face equivalent source approximation. In the approximations of $\mathcal{A}_{B}(v)(\bsx)$ by high-order quadratures, contributions coming from source points that lie outside  $\mathcal{N}_{\bsx}$ is obtained accurately by replacing these sources with equivalent sources on a sparse regular grid for efficient computations using FFT. Note that the source locations include points in $\Omega_{B,h}$ (depicted as square shape ) as well as those coming from $\Omega_{h} \setminus \Omega_{B}$ (depicted as cross shape) that do not contribute to the boundary integral. The source strength corresponding to such virtual sources are, therefore, set to zero. This is discussed in more detail in Section \ref{sec:-nsnai}}.
		 \label{fig:-accelpic}
	\end{center}

\end{figure}

For  $\bsx_q \in \Omega_{h}$, the non-singular integrals in (\ref{eq:-BrdySplit}) can be  approximated to high-order by means of any classical high-order quadrature, for example, the trapezoidal rule for $s'$ integration and composite Newton-Cotes in $t'$ direction that we utilized earlier. This application of quadrature, obviously, have contributions coming from all of $\Omega_{B,h}$ including those in $\mathcal{N}_{\bsx_q}$ that we call {\em adjacent contributions} and those from outside that we refer to as {\em non-adjacent contributions}. The non-adjacent interaction, in particular, takes a simple form of a discrete convolution that reads
\begin{equation}
\label{eq:-nans}
\mcA_{na,B}^{reg} (v)(\bsx_q) = \sum_{\bm{y}_{l}\in\Omega_{B,h}\setminus \mathcal{N}_{\bsx_{q}}} w_lG_{\kappa}(\bsx_q - \bm{y}_l) v(\bm{y}_l) \eta_{\tau_0}((\mcP\bm{y}_l-\bm{y}_l)\cdot \bsnu(\mcP\bm{y}_l)).
\end{equation}
While a discrete convolution can be computed efficiently using FFT if the underlying computational grid is regular, we note that, in (\ref{eq:-nans}), the contributing sources, that we refer as ``true sources", are placed somewhat irregularly in $\Omega_{B}$. In addition, we also note that while the sources range in $\Omega_{B,h}$, the convolution needs to be evaluated not just on points in $\Omega_{B,h}$, but also on other grid points coming from $\Omega_h \setminus\Omega_{B}$ that we refer to as ``virtual sources". These additional sources, obviously, do not contribute to the sum in (\ref{eq:-nans}) but are included in the computation with ``zero" strength (weight) so that the discrete convolutions get evaluated at the target points $x_q \in \Omega_h \setminus\Omega_{B}$ in an efficient manner using our FFT based acceleration strategy \cite{Anand2015highorder,bruno2001fast} that we explain next for completeness.

The acceleration strategy seeks to compute a discrete convolution 
by substituting  true sources by a certain set of ``equivalent sources" located at a regular Cartesian grid on the two parallel faces of cell $c_{ij}$ (see Figure \ref{fig:-accelpic}). 
Let $\sigma_{{ij,l}}^{m}$, $\sigma_{{ij,l}}^{d}$ $\left(l=1,\cdots, N^{eq}\right)$ denotes acoustic monopoles and dipoles respectively, at the Cartesian grid $\bsx_{ij,l}$ on the parallel faces of cell $c_{ij}$. Then the field generated by these equivalent sources is given by
\begin{equation}
\label{eq:-fges}
\mcA_{{ij},B}^{reg,eq} (v)(\bsx) = \sum_{l=1}^{N^{eq}} \left(\sigma_{ij,l}^{m}G_{\kappa}(\bsx - \bsx_{ij,l})+ \sigma_{ij,l}^{d}\frac{\partial G_{\kappa}(\bsx - \bsx_{ij,l})}{\partial \bm{\nu}\left( \bsx_{ij,l}\right)}\right).
\end{equation}
Let $\mcA_{ij,B}^{reg,true} (v)(\bsx)$ denote the field generated by true sources within cell $c_{ij}$, say, given by
\begin{equation}
\label{eq:-fgts}
\mcA_{ij,B}^{reg,true} (v)(\bsx) = \sum_{\bm{y}_{l}\in c_{ij} }w_lG_{\kappa}(\bsx - \bm{y}_l) v(\bm{y}_l) \eta_{\tau_0}((\mcP\bm{y}_l-\bm{y}_l)\cdot \bsnu(\mcP\bm{y}_l)).
\end{equation}
The unknown quantities $\sigma_{ij,l}^{m}$ and  $\sigma_{ij,l}^{d}$ are obtained as a solution of over determinant linear system $\bm{A\sigma}=\bm{b}  $. Here, the matrix $\bm{A}$ is obtained by evaluating (\ref{eq:-fges}) at $n^{coll} =4N^{eq}$ points located on the boundary of adjacent set $\mathcal{N}_{\bsx}$ and the vector $\bm{b}$ correspond to the evaluation of (\ref{eq:-fgts}) at $n^{coll}$ points. It is important to note that,  because of the identical geometry of every cell $c_{ij}$, the $\bm{QR}$ decomposition of matrix $\bm{A}$ is obtained once and stored for repeated use.  This whole process require $O(N^{3/2}/L^{3}) +O(N^{3/2}/L^2) $ operations in total.

It is known  that, as long as $\bsx$ is non-adjacent to cell $c_{ij}$,  $\mcA_{ij,B}^{reg,eq} (v)(\bsx)$ provide an accurate approximation to $\mcA_{ij,B}^{reg,true} (v)(\bsx)$ (see \cite{bruno2001fast}).  To be more precise, for any $\bsx \in c_{ij}$, the quantity
\begin{equation}
\label{eq:-eqsc}
\mcA_{B}^{reg,eq} (v)(\bsx) = \sum_{k=1}^{L} \sum_{l=1}^{L}\mcA_{{kl},B}^{reg,eq} (v)(\bsx)- \sum_{k=i-1}^{i+1} \sum_{l=j-1}^{j+1}\mcA_{{kl},B}^{reg,eq} (v)(\bsx)
\end{equation}
  provide an accurate approximation of non-adjacent interaction $\mcA_{na,B}^{reg} (v)(\bsx)$ defined in (\ref{eq:-nans}). The importance of (\ref{eq:-eqsc}) lies in the fact that, it is a 
  convolution over a regular Cartesian grid, and hence,  for all equivalent source points, can be evaluated in a fast way by means of FFT in $O\left(LN^{1/2}\log(LN^{1/2})\right)$ operations.
  
  Finally, to evaluate field values at true source locations or at any virtual source point, a free space Helmholtz equation within each cell $c_{ij}$ with Dirichlet boundary data coming from  $\mcA_{B}^{reg,eq} (v)(\bsx) $ is solved. To ensure unique solvability of the problems, we assume that the cells are non-resonant. To obtain the solution of these well posed problems efficiently, a discretized plane wave expansion in a coordinate system  local to cell $c_{ij}$, of the form  
  \begin{equation}
  \label{eq:-pwe}
  \mcA_{B}^{reg,eq} (v)(\bsx) \approx \sum_{l=1}^{N^{coll}} \beta_{l}\exp \left(i k\bm{d}_{l}(\bsx-\bsx_{ij,c}) \right)
  \end{equation}  
  is used. Here, $\bm{d}_{l}$ denote the unit vector on the surface of disc of unit radius centered at $\bsx_{ij,c}$, center of the cell $c_{ij}$. The coefficients $\beta_l$ are obtained by solving an overdetermined linear system which is constructed by enforcing (\ref{eq:-pwe}) at equivalent source points located on the boundary of cell $c_{ij}$. Note that, series (\ref{eq:-pwe}) converges spectrally provided $N^{coll} =O(N^{1/2}/L)$, as discussed in \cite{Anand2015highorder,bruno2001fast}.
  Again, owing to identical geometry of cell $c_{ij}s$, the $\bm{QR}$ factorization of coefficient matrix needs to be computed only once and stored for repeated use. Thus, evaluation of $ \mcA_{B}^{reg,eq} (v)(\bsx)$ at all true and virtual source locations could be completed in an $O\left(N^{3/2}/L^{3}\right)+O\left(N^{3/2}/L\right)$ operations. Thus, the overall computational cost for non-adjacent interactions stands at $O(LN^{1/2}\log (LN^{1/2})) +O\left(N^{3/2}/L^{3}\right)+O\left(N^{3/2}/L\right)$.

\section{Boundary-base interaction}
\label{sec:-BBulkI}
Indeed, the algorithm proposed for base-base interactions, in section \ref{sec:-BubI}, computes the base integral operator $\mcA_E (v)(\bsx) $ only for $\bsx \in \Omega_{E,h}$. This calculation for other points, obviously, can be accomplished by in a straightforward manner by directly computing the integration 
in $\mcA_E (v)(\bsx) $. We, however, avoid this compatitively expensive direct computation by obtaining these values through an interpolation of the base-base data. While the exact details of interpolation strategy, of course, depend on the underlying base grid $\Omega_{E,h}$, it relies of the following two requirements: 
first, $\Omega_{E}$ contains $\Omega_{B}$, and second, smoothness of the 
volume potential $\mcA_E (v)(\bsx)$. 
For example,
%
%
%
in the Addition theorem approach, a piece-wise Chebyshev radial gird and a uniform grid in angular direction constitutes the computational grid, whereas, the pre-corrected  trapezoidal  approach utilizes a regular grid in all Cartesian directions. While the key idea underlying our interpolation scheme is common to both these cases, and is closely related to ideas presented in \cite{bruno2001fast}, the finer details do differ due to differences in the grid structure. We highlight both these distinct strategies in what follows starting with details in the pre-corrected trapezoidal approach. 

Choose a smooth periodic window function $\omega(\bsx)$ such that it takes value one in $\Omega_R$ and varies smoothly to zero towards the boundary of $\Omega_{E}$, where  the set $\Omega_R$  is rectangular domain such that $\Omega \subset \Omega_{R} \subset \Omega_{E}$. As $\omega(\bsx)\mcA_E (v)(\bsx) =\mcA_E (v)(\bsx),$ $\forall \bsx \in \Omega $, we can approximate $\mcA_E (v)(\bsx)$ at any arbitrary point $\bsx \in \Omega$ by interpolating $\omega(\bsx)\mcA_E (v)(\bsx)$ using precomputed data $\mcA_E (v)(\bsx), \bsx \in \Omega_{E,h}$.  Now, to interpolate smooth periodic function $\omega(\bsx)\mcA_E (v)(\bsx)$, we use FFT-Refined polynomial interpolation strategy \cite{Anand2015highorder} what we summmarize as the followowing three step procedure: 

\begin{enumerate}
	\item   Obtain the Fourier coefficient of   $\omega(\bsx)\mcA_E (v)(\bsx)$ through FFT using its value on a uniform grid.
	\item  Evaluate the Fourier series of  $\omega(\bsx)\mcA_E (v)(\bsx)$ on a refined regular grid. 
	\item Construct local interpolating polynomials of  a fixed degree using values of  $\omega(\bsx)\mcA_E (v)(\bsx)$ on the refined grid obtained in Step 2. 
\end{enumerate}

For the Addition theorem approach, on the other hand, 
where the base integral 
operator data for $\mcA_E (v)\left(\bsx(r,\theta)\right)$ is available on a regular $\theta$-grid and a piecewise Chebyshev $r$-grid,
the refined grid data
is obtained via Fourier and Chebyshev interpolants in $\theta$ and $r$ respectively. As in the pre-corrected trapezoidal case, local polynomial interpolation is then used for evaluation at arbirary points in $\Omega_B$.
\section{Base-boundary interaction}
\label{ref:-bulk-boundary}
In this section, we elucidate our approach for computation of the integral $\mcA_B (v)(\bsx)$, defined in (\ref{eq:-ABdry}), for the case when 
$\bsx \in \Omega_{E,h}$. 
For ease of explanation, we partition $\Omega_{E,h}$ into three disjoint sets, first $\Omega_{E,h}^{e}$ contains points in the exterior of $\Omega$ ( $\Omega_{E,h}^{e} = \Omega_{E,h} \cap \Omega^c$), second $\Omega_{E,h}^{b}$ contains boundary points ($\Omega_{E,h}^{b} = \Omega_{E,h} \cap \Omega_B$) and the last $\Omega_{E,h}^{i} = \Omega_{E,h} \cap \Omega \cap \Omega_{B}^c$ contains the rest.

As we have pointed out in subsection \ref{sec:-Dis_int}, the approximation grid $\Omega_{h}$ is contained in $\Omega$, and therefore, $\Omega_{h}$ does not overlap with $ \Omega_{E,h}^{e}$. 
When the target point $\bsx \in \Omega_{E,h}^{b}  $, it typically does not  coincide with any point in $\Omega_{B,h}$, and a direct application of integration scheme employed for boundary-boundary interaction, discussed in section \ref{sec:-BBI}, is neither accurate nor efficient.  To address this, we use a strategy similar to the one discussed in section \ref{sec:-BBulkI}.
More precisely, for any off grid target point $\bsx \in \Omega_{B}\setminus \Omega_{B,h}$,  we obtain high-order approximation of $\mcA_B (v)(\bsx)$ by interpolating its pre-computed values on the boundary grid $\Omega_{B,h}$.    However, as $\mcA_B (v)(\bsx)$ is not compactly supported, 
we utilize partition of unity again toward obtaining high order approximations. Using the partition of unity $\left\{w_k :k = 1,...,K\right\}$ subordinate to the covering $\{ \Omega_{k,B} \}_{k = 1}^K$ of $\Omega_B$, we write
\begin{equation}
\mcA_B (v)(\bsx)=\sum_{k \in \mathcal{I}_{B}\left(\bsx \right)}\left(w_k\mcA_B (v)\right)(\bsx),
\end{equation}
  where the index set is given by $I_{B}(\bsx)=\left\{k\ \vert\ \bsx \in \Omega_{k,B} \right\}$. As
   $\left(w_k\mcA_B (v)\right)(\bsx)= \left(w_k\mcA_B (v)\right)\left(\bsxi_k(s,t)\right)$ is
   periodic in $s$-variable and smooth in both the variable $s$ and $t$,  we can approximate  $\left(w_k\mcA_B (v)\right)(\bsx)$ to high-order accuracy by means of  interpolation technique similar to the one discussed in section \ref{sec:-BBulkI} for boundary-base interactions. 
   We note that, as boundary region $\Omega_{B}$ is thin in the transverse direction, trigonometric refinement is necessary only in $s$-variable.
   Again, we evaluate  $ \left(w_k\mcA_B (v)\right)\left(\bsxi_k(s,t)\right)$ on finer $s$-grid lines by means of FFT as $ \left(w_k\mcA_B (v)\right)\left(\bsxi_k(s,t)\right)$
   is smooth and periodic in $s$-variable and, subsequently, approximate $\left(w_k\mcA_B (v)\right)\left(\bsxi_k(s,t)\right)$ locally to high-order by interpolating polynomials of a fixed degree based on the refined grid data. Due to relatively small size of $\Omega_{E,h}^b$, this interpolation procedures contribution toward the overall computational time is insignificant when compared with other components of the numerical scheme.

 Finally, for the target point $\bsx \in \Omega_{E,h}^{i}$, as the integral  $\mcA_B (v)(\bsx)=\mcA_B \left(v)\right)\left(\bsxi_k(s,t)\right)$ has smooth integrand, but with periodicity only in $s-$variable, a
   high-order approximation is obtained by employing trapezoidal rule in $s-$variable and composite Newton-Cotes in $t-$variable.  This straight forward application of classical quadratures results in $O(NN_{B})$ computational cost, where $N_B$ is the number of unknowns coming from the discretization of boundary region $\Omega_{B}$.
   This computation, however, could be expedited further by means of two face equivalent source approximation technique that we discussed in section \ref{sec:-BBI}. Toward this, we divide evaluation of  $\mcA_B (v)(\bsx)$ into adjacent and non-adjacent interactions where former is obtained accurately by means of a high-order quadrature while the later is approximated accurately and efficiently through the plane wave expansion formula (\ref{eq:-pwe}). The computational cost, therefore, is determined by the cost of evaluation of formula (\ref{eq:-pwe}) and the cost of computing adjacent interactions for all target points in $\Omega_{E,h}^{b}$ and  $\Omega_{E,h}^{i}$.  As each cell $c_{ij}$ contains, on an average, $O(N_{B}/L^2)$ boundary discretization points, the adjacent calculation 
   requires only $O(NN_{B}/L^2)$ operations. The cost of evaluation of non-adjacent contributions, on the other hand, with the choice $N^{coll} =O(N^{1/2}/L)$ in (\ref{eq:-pwe}) stand at $O(N^{3/2}/
   L)$.
%
%

\section{Computational Cost}
\label{sec:-Comp-Cost}
The computational cost of methods considered in the proposed framework is determined by the cost of evaluation of volume potential
 $\mathcal{A}(v)(\bsx)$ in (\ref{eq:-vol-pot}), at all grid points in $\Omega_{h}$. As we have discussed in section \ref{sec:-OS}, evaluation of integral operator (\ref{eq:-vol-pot}) at all grid points are completed by evaluating four different interactions. 	An unaccelerated computation of boundary-boundary and base-boundary interactions would result in an over all computational complexity  $O(NN_{B})$. However, as the thickness of the boundary region $\Omega_{B}$ is of the order of wavelength, we can safely take $N_{B}=O(N^{1/2})$, resulting in the complexity $O(N^{3/2})$.  This could be further improved to $O(N \log N)$ cost by accelerating boundary-boundary and base-boundary interactions using techniques described in section \ref{sec:-BBI}. In the following, we briefly provide a complexity analysis to elaborate on this claim:  
\begin{itemize}
\item As we have elucidated in section \ref{sec:-BubI}, base-base interaction using either approach,  Addition theorem or pre-corrected trapezoidal, is obtained by means of FFT and, therefore, requires only $O(N \log N)$ operations. 
\item Evaluation of boundary-boundary interaction is completed with a total computational cost of  $O(LN^{1/2}\log (LN^{1/2})) +O\left(N^{3/2}/L^{3}\right)+O\left(N^{3/2}/L\right)$. By choosing, parameter $L=O(N^{1/2})$, computational  cost of this step reduces to $O(N \log N)$. A detailed discussion on this can be found in \cite{Anand2015highorder}.
\item Use  of FFT-refined polynomial interpolation scheme for evaluation of base-boundary interaction results in  $O(N\log N)$ computational cost.  

\item As discussed toward the end of section \ref{ref:-bulk-boundary}, the base-boundary interaction requires $O\left(NN_{B}/L^2\right) +O\left(N^{3/2}/L\right)$ operations. With the choice $L =O(N^{1/2})$, that we make, the computational cost for this component of the algorithm exhibits $O(N )$ computational complexity.

    \end{itemize}  

Summing up the costs arising out of all four interactions, we conclude that the total computational complexity of the methods falling within this framework  stand at $O(N\log N)$.

\section{Numerical Results}
\label{sec:-Num-Res}
In this section, we demonstrate the enhanced high-order accuracy of the two methods considered under the proposed framework  through a variety of computational examples.
The numerical results presented in this section  are obtained using C$++$ implementations of our schemes. The relative error (in the near field) reported here are  computed as  
\begin{align*}
\varepsilon_{\infty} &= \frac{ \underset{1\leq i \leq N}{\max} \left|u^{\text{exact}}(\bm{x}_{i})-u^{\text{approx}}(\bm{x}_{i})\right|}{ \underset{1\leq i \leq N}{\max} \left|u^{\text{exact}}(\bm{x}_{i}) \right|}, \\
\varepsilon_{2} &= \left ( \frac{\sum \limits_{i=1}^{N}\left|u^{\text{exact}}(\bm{x}_{i})-u^{\text{approx}}(\bm{x}_{i})\right|^2}{\sum \limits_{i=1}^{N}\left|u^{\text{exact}}(\bm{x}_{i}) \right|^2} \right)^\frac{1}{2}.
\end{align*}
 We use the notation $P_{1}\times N_{1}\times N_{2}+ M_1\times M_2$ to specify that $P_1$ number of overlapping boundary patches, each with $N_{1}\times N_{2}$ discretization points and $M_{1}\times M_{2}$ points over base region are used for the corresponding numerical solution.
 In all the tabulated results, the acronym numIt. denotes the number of GMRES iterations required to achieve the desired accuracy and ``Order" denotes the numerical order of convergence.
 
  As discussed in the introduction, this manuscript provides a technique to improve the order of convergence of those methods that converge to high-order only for smoothly varying media and result in low order accuracy for discontinuous scattering media. For specific examples of this, we have considered two algorithms, one based on the Addition theorem approach \cite{bruno2004efficient}, and the other utilizing the  pre-corrected trapezoidal rule  \cite{aguilar2004high,duan2009high}, both of which indeed converge with high-order for smooth scatterers, while yielding low order accuracy for discontinuous scattering media.
  In view of this,  all scattering calculations in this section, correspond to only discontinuous media. 
  We note that, in all examples, $5-$points composite Newton-cotes quadrature is used  
   to approximate the transverse integral over the boundary region in (\ref{eq:-BdrySplitK}). 
   
  
\begin{exmp}(\textit{A convergence study for proposed integration scheme})
\end{exmp}
As we have explained,  the proposed high-order method  relies on  high-order evaluation of volume potential $\mathcal{A}(v)$ defined in (\ref{eq:-vol-pot}). Therefore,  in our first example, we present numerical results to corroborate the high-order convergence of the two methods implemented under this framework. 

We consider a disc of acoustical size $\kappa a=4$ ($a$ is the diameter of disc ),  with the refractive index   $n(\bsx) =\sqrt{2}$ when $\bsx$ happens to be within the disc and one otherwise. In this case, using Addition theorem for Hankel kernel and plane wave expansion for incident wave \cite{colton2012inverse},    
volume potential  $\mathcal{A}_{E}(v)$ expressed analytically. In order to show the convergence, error incurred in our numerical approximation at different level of discretization are reported in Tables \ref{table:-PCTFMK2} and  \ref{table:-ADTHK2}. The results in Table \ref{table:-PCTFMK2} and  \ref{table:-ADTHK2}  are obtained by employing pre-corrected trapezoidal rule and Addition theorem method under the proposed framework respectively. These studies clearly illustrate the enhanced rate of convergence for the two methods under investigation when augmented as we proposed in this paper.  


\begin{table}  [h!]
	\begin{center} 
		\begin{tabular}{c|c|c|c|c|c}
			\hline
			Grid Size & Unknowns  & \multicolumn{2}{c}{$L^{\infty}$}&\multicolumn{2}{c} {$L^2$} \\ 
			\cline{3-6} 
			&   & $\varepsilon_{\infty}$ & Order&  $\varepsilon_{2}$ & Order \\ 
			\hline
			$2\times9\times5+17\times17 $&379&4.67e-02&-&3.73e-02& - \\ 
			$2\times17\times9 +33\times33$&1395&3.00e-03&3.96e+00&2.48e-03&3.91e+00\\ 
			$2\times33\times17 +65\times65$&5347&3.91e-04&2.94e+00&3.40e-04&2.86e+00\\ 
			$2\times65\times33 +129\times129$&20931&1.52e-05&4.68e+00&7.72e-06&5.46e+00\\ 
			$2\times129\times65 +257\times257$&82819&3.91e-07&5.28e+00&1.30e-07&5.90e+00\\
			\hline 
		\end{tabular}  
	\end{center} 
	\caption{Convergence study: Approximation of volume potential $\mathcal{A}(v)$
		using pre-corrected trapezoidal rule under the proposed framework for a
		 disc shape scatterer with  $\kappa a=4$ and $n=\sqrt{2}$. } 
		\label{table:-PCTFMK2}
\end{table} 

\begin{table} [h!]

	\begin{center} 
		\begin{tabular}{c|c|c|c|c|c}\hline
			
			Grid Size & Unknowns  & \multicolumn{2}{c|}{$L^{\infty}$}&\multicolumn{2}{c} {$L^2$} \\ 
			\cline{3-6}&   & $\varepsilon_{\infty}$ & Order&  $\varepsilon_{2}$ & Order \\ 
			\hline
			$2\times9\times5+17\times 9 $&243&1.25e-01&-&1.08e-01& - \\ 
			$2\times17\times9 +33\times 17$&867&3.61e-03&5.12e+00&1.24e-03&6.45e+00\\ 
			$2\times33\times17 +65\times 33$&3267&2.65e-04&3.77e+00&8.02e-05&3.95e+00\\ 
			$2\times65\times33 +129\times 65$&12675&6.65e-06&5.32e+00&5.00e-06&4.00e+00\\ 
			$2\times129\times65 +257\times 129$&49923&7.08e-08&6.55e+00&5.43e-08&6.52e+00\\ 
			\hline
		\end{tabular}  
	\end{center} 
	\caption{Convergence study: Approximation of volume potential $\mathcal{A}(v)$
		using Addition theorem method under the proposed framework for a
		disc shape scatterer with $\kappa a=4$ and $n=\sqrt{2}$. } 
		\label{table:-ADTHK2}
\end{table} 
\begin{table}[h!]	
	\begin{center}
		\begin{tabular}{c|c|c|c|c|c|c} \hline
			$\kappa a$ & \multicolumn{3}{c|}{pre-corrected trapezoidal rule} & \multicolumn{3}{c}{current}  \\ 
			\hline	
			\cline{2-7}& $N$ & $\varepsilon_{\infty}$  & $\varepsilon_{2}$  & $N$ & $\varepsilon_{\infty}$  & $\varepsilon_{2}$ \\ 			
			\hline
			20 & $4198401$ & 2.68e-05 & 2.10e-04  & 74563& 7.29e-07&6.26e-07 \\		
			\hline
			30 & $4198401$ & 3.53e-04 & 2.55e-04  &74563 &1.60e-06&1.19 e-06  \\
			\hline
			%
		\end{tabular}
		\caption{ Accuracy in the approximation of volume potential $\mathcal{A}(v)$ computed by original pre-corrected trapezoidal rule and the current algorithm.} 
		\label{table:-PCTVC}
	\end{center}	
\end{table}	

%
In Table \ref{table:-PCTVC}, we compare the accuracy of our quadrature scheme with the original pre-corrected trapezoidal rule. We see that, the present algorithm achieves considerably higher accuracy than the pre-corrected trapezoidal rule with relatively fewer unknowns. Moreover, the accuracy gap widens rapidly as the number of unknowns employed for the approximation increases. For instance, in the calculations with $\kappa a=30$, the current algorithm produces more than $220$ times  better accuracy with $56$ time fewer unknowns when compare to what results from the pre-corrected trapezoidal rule . 

\begin{exmp}(\textit{Convergence study for a simple scatterer})
\end{exmp}
As a second exercise,  using our high-order integration scheme, we simulate the scattering by a penetrable disc of acoustical size $\kappa a=4\pi$ ($a$ being diameter of the inhomogeneity) for which the true solution can be evaluated analytically.
We again take $ n(\bsx) =\sqrt{2}$ if  $\bm{x} \in \Omega$ and one otherwise.
Note that $n(\bsx)$ is discontinuous across the interface of scatterer $\Omega$.  
 To obtain the numerical solution of this problem, boundary region is covered by two thin overlapping annular  patches. In order to demonstrate  rapid convergence of our algorithm, numerical  solutions are computed on several levels of discretization and the corresponding results are tabulated in Tables \ref{tabel:-PCTK2PIDISC} and \ref{tabel:-ADTHK2PIDISC}.

In Table \ref{tabel:-PCTK2PIDISC}, we present  numerical results  corresponding to the  pre-corrected trapezoidal rule under the proposed  framework, while those coming from Addition theorem approach are given in Table \ref{tabel:-ADTHK2PIDISC}. 

\begin{figure}[!b] \label{fig:-Disc-K25}
	\begin{center}	
		\subfigure[Absolute value of the scattered field, $|\mathfrak{u}^s|$]	{\includegraphics[clip=true, trim=000 00  00 00, scale=0.20]{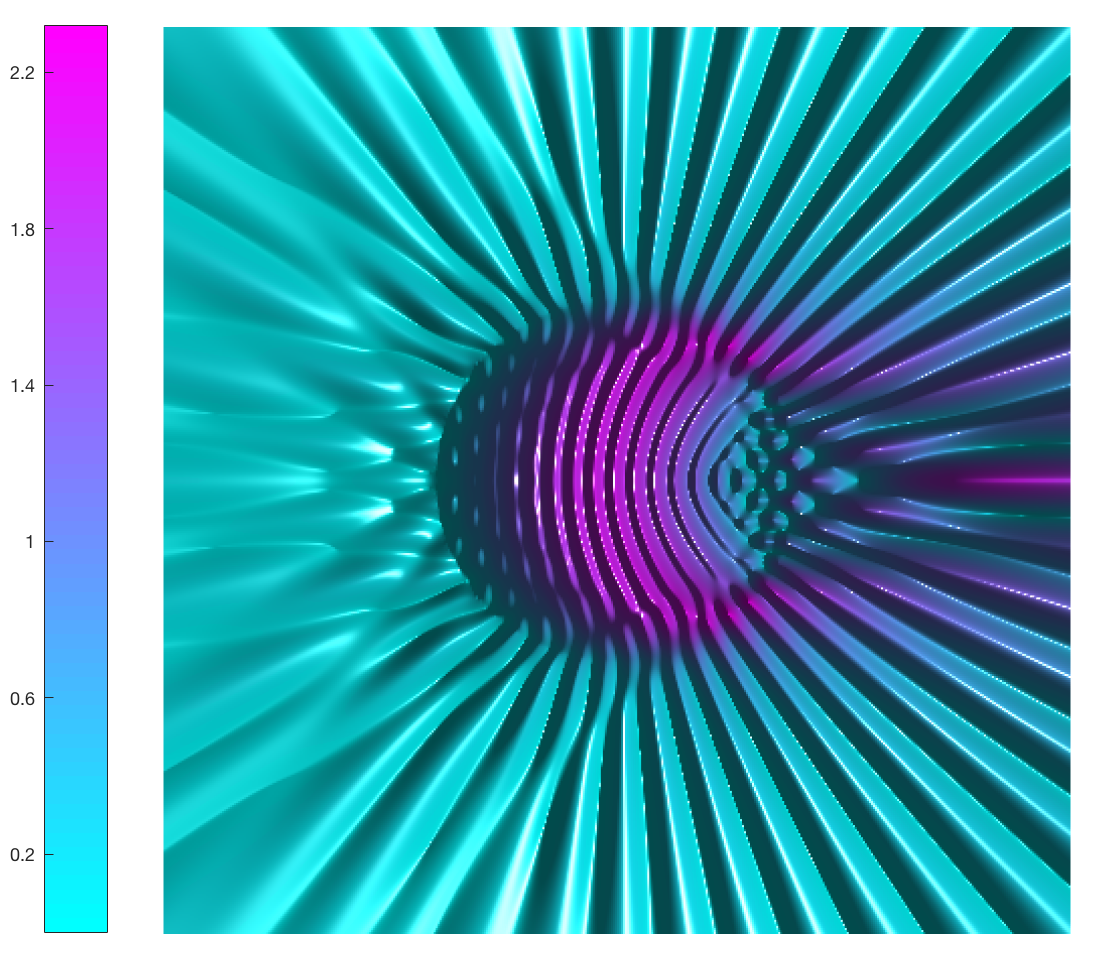}}
		\hfill
		\subfigure[Absolute value of the total field, $|\mathfrak{u}|$]{\includegraphics[clip=true, trim=000 0  0 00, scale=0.20]{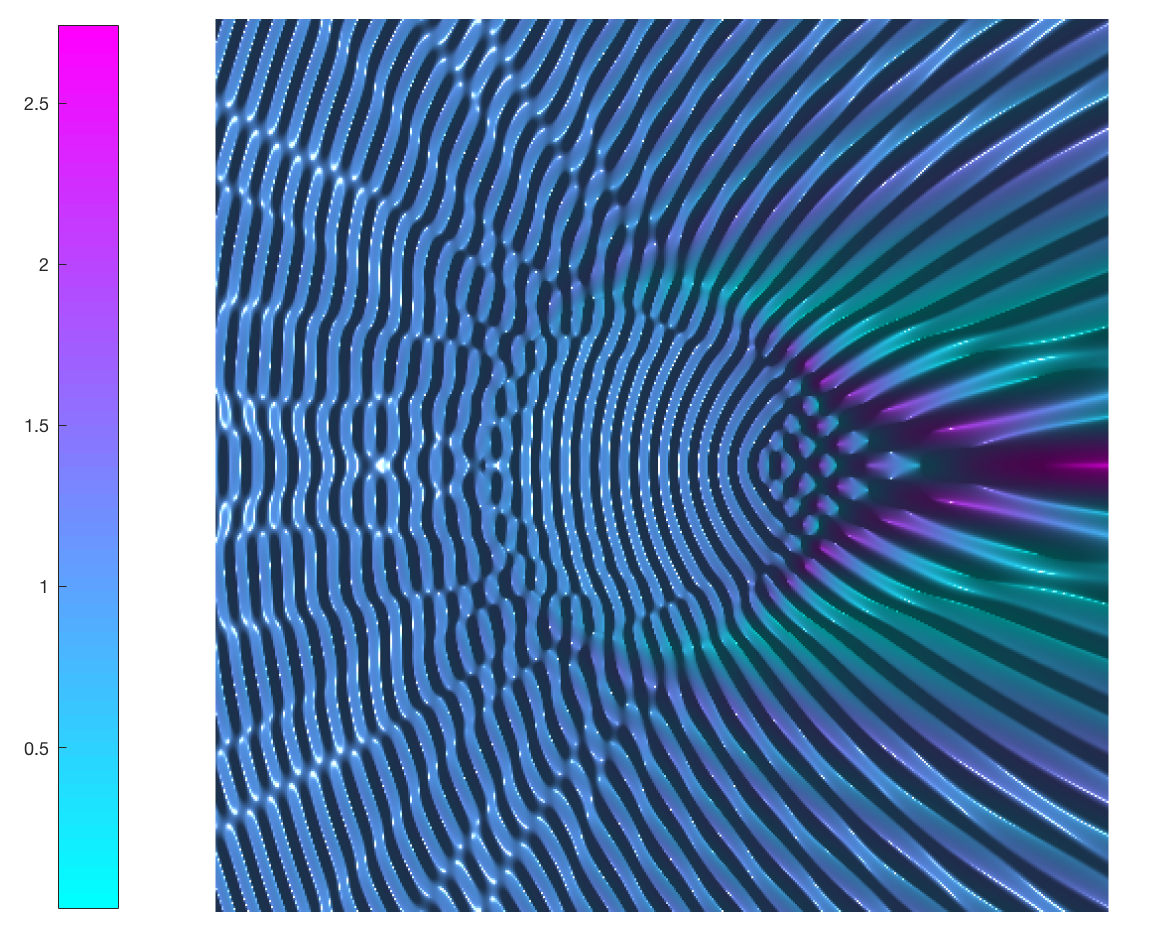}}
		\hfill
		\caption{Scattering of a plane wave $ \exp (i \kappa x)$ by a penetrable  disc of acoustical size $\kappa a =50$.
			Fields are computed using Addition theorem method under the proposed framework on the grid of size  $2\times{65} \times{129}+  129 \times 65 $. In the near field  we obtain an accuracy of 0.001.}
	\end{center}
\end{figure}

For a pictorial visualization, in Figure \ref{fig:-Disc-K25}, we plot absolute value of scattered and total fields for the penetrable disc of acoustical size $\kappa a =50$ with contrast function $m(\bsx) =-0.3$ inside the disc and zero elsewhere.  For this experiment the numerical solution is obtained using our approach with  Addition theorem where we have three digits of accuracy.
 
\begin{table}[!h]
	\begin{center} 
	\begin{tabular}{c|c|c|c|c|c|c} \hline
		Grid Size & Unknowns  & \multicolumn{2}{c|}{$L^{\infty}$}&\multicolumn{2}{c|} {$L^2$} &  numIt\\ 
		\cline{3-6}&   & $\varepsilon_{\infty}$ & Order&  $\varepsilon_{2}$ & Order&  \\ 
		\hline
		$2\times9\times5+17\times17 $&379&7.29e-01&-&5.74e-01& - &8\\ 
		$2\times17\times9 +33\times33$&1395&7.88e-02&3.21e+00&9.23e-02&2.64e+00&17\\ 
		$2\times33\times17 +65\times65$&5347&5.92e-03&3.74e+00&5.60e-03&4.04e+00&22\\ 
		$2\times65\times33 +129\times129$&20931&1.46e-04&5.35e+00&1.68e-04&5.06e+00&30\\ 
		$2\times129\times65 +257\times257$&82819&2.51e-06&5.86e+00&2.89e-06&5.86e+00&31\\ 
		\hline
	\end{tabular}  
	\caption{Convergence study: Plane wave scattering by a penetrable disc with  $\kappa a =4 \pi$ and refractive index $n(\bm{x}) = \sqrt{2}$ when $\bm{x} \in \Omega$ and one otherwise, where the base
	integral operator is approximated by means of pre-corrected trapezoidal rule.} 
	\label{tabel:-PCTK2PIDISC} 
\end{center} 
\end{table} 

\begin{table} [!h]
	\begin{center} 
		\begin{tabular}{c|c|c|c|c|c|c}\hline
			Grid Size & Unknowns  & \multicolumn{2}{c|}{$L^{\infty}$}&\multicolumn{2}{c|} {$L^2$} &  numIt\\ 
			\cline{3-6}
			&   & $\varepsilon_{\infty}$ & Order&  $\varepsilon_{2}$ & Order&  \\ 
			\hline
			$2\times9\times5+17\times9 $&243&4.54e-01&-&3.56e-01& - &8\\ 
			$2\times17\times9 +33\times17$&867&3.73e-02&3.61e+00&4.85e-02&2.87e+00&16\\ 
			$2\times33\times17 +65\times33$&3267&2.70e-03&3.79e+00&2.34e-03&4.37e+00&25\\ 
			$2\times65\times33 +129\times65$&12675&1.70e-04&3.99e+00&7.34e-05&5.00e+00&30\\ 
			$2\times129\times 65+257\times129$&49923&2.15e-06&6.31e+00&4.98e-07&7.20e+00&35\\ 
			\hline
		\end{tabular}  
		\caption{Convergence study: Plane wave scattering by a penetrable disc with $\kappa a =4 \pi$ and refractive index $n(\bm{x}) = \sqrt{2}$ when $\bm{x} \in \Omega$ and one otherwise, where
		base integral operator is approximated by means of	
			Addition theorem method.}
		 \label{tabel:-ADTHK2PIDISC}
	\end{center} 
\end{table}

\begin{exmp}(\textit{Convergence study for a complex scatterer})
\end{exmp}
In the previous example, we have considered a disc shape scatterer with 
constant inhomogeneity. However,  our algorithm is not restricted to either the simplicity of scattering geometry or  to  constant material properties.  
 This example will demonstrate the adaptability and applicability of the proposed method in dealing with scatterers that have relatively complex geometrical description as well as variable martial properties.
Towards this, we consider scattering by penetrable inhomogeneous bean shaped scatterer, as depicted in Figure \ref{fig:-omega_spliting}, whose boundary curve is defined by $
r(t) = \left(\cos t + 0.65 \cos 2t -0.65, 1.5 \sin t \right), \ \ t \in [0,2\pi). $  

For the numerical approximation, we cover the boundary region by two overlapping patches, as shown in Figure \ref{fig:-boundarypou}. 
To study the convergence behavior, once again, we compute the total field $u$ under the plane wave incidence $u^{i} = \exp(i\kappa x)$. As an analytical solution is not available for this scattering configuration, we use numerical solution obtained by our algorithm on a finer grids for comparisons. We take acoustical size of the scatterer  $\kappa a =10 \pi$, and contrast function $m(\bsx)$ is given by
\begin{equation}\label{eq:-refI}
 m(\bm{x}) =
\begin{cases}
1-0.5\exp\left(-\left(x^2+y^2\right)\right)& \mbox{if} \ \ \bm{x}  = (x,y) \hspace{1mm}\in \Omega \\
0 & \mbox{if} \ \ \bm{x}=(x,y) \in \mathbb{R}^{2}\setminus\Omega. \\
\end{cases}
\end{equation}
Note that, the contrast function $m(\bsx)$ is discontinuous across the scattering interface $\partial\Omega$. 
The numerical results corresponding to pre-corrected trapezoidal rule and Addition theorem method are presented in Table \ref{table:-PCTKA10PIBEAN} and \ref{tabel:-ADTHKA10PIBEAN} respectively. 
These computational results confirm that the high-order convergence of our method remain valid for complicated geometries with discontinuous variable material properties.

\begin{figure}[h!] \label{fig:-bean-K25}
	\begin{center}	
		\subfigure[Absolute value of the scattered field, $|\mathfrak{u}^s|$]	{\includegraphics[clip=true, trim=000 00  00 00, scale=0.19]{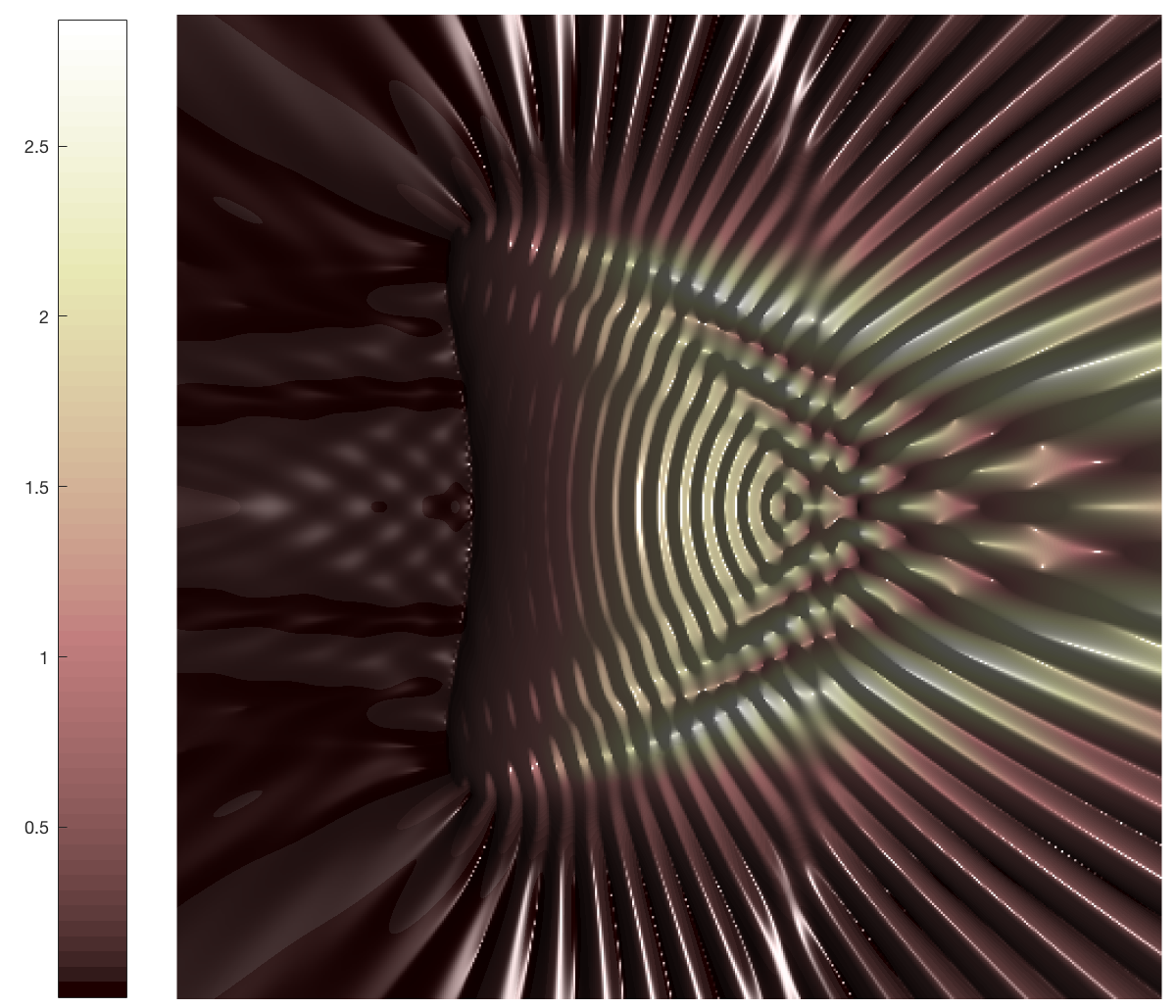}}
		\hfill
		\subfigure[Absolute value of the total field, $|\mathfrak{u}|$]{\includegraphics[clip=true, trim=000 0  0 00, scale=0.19]{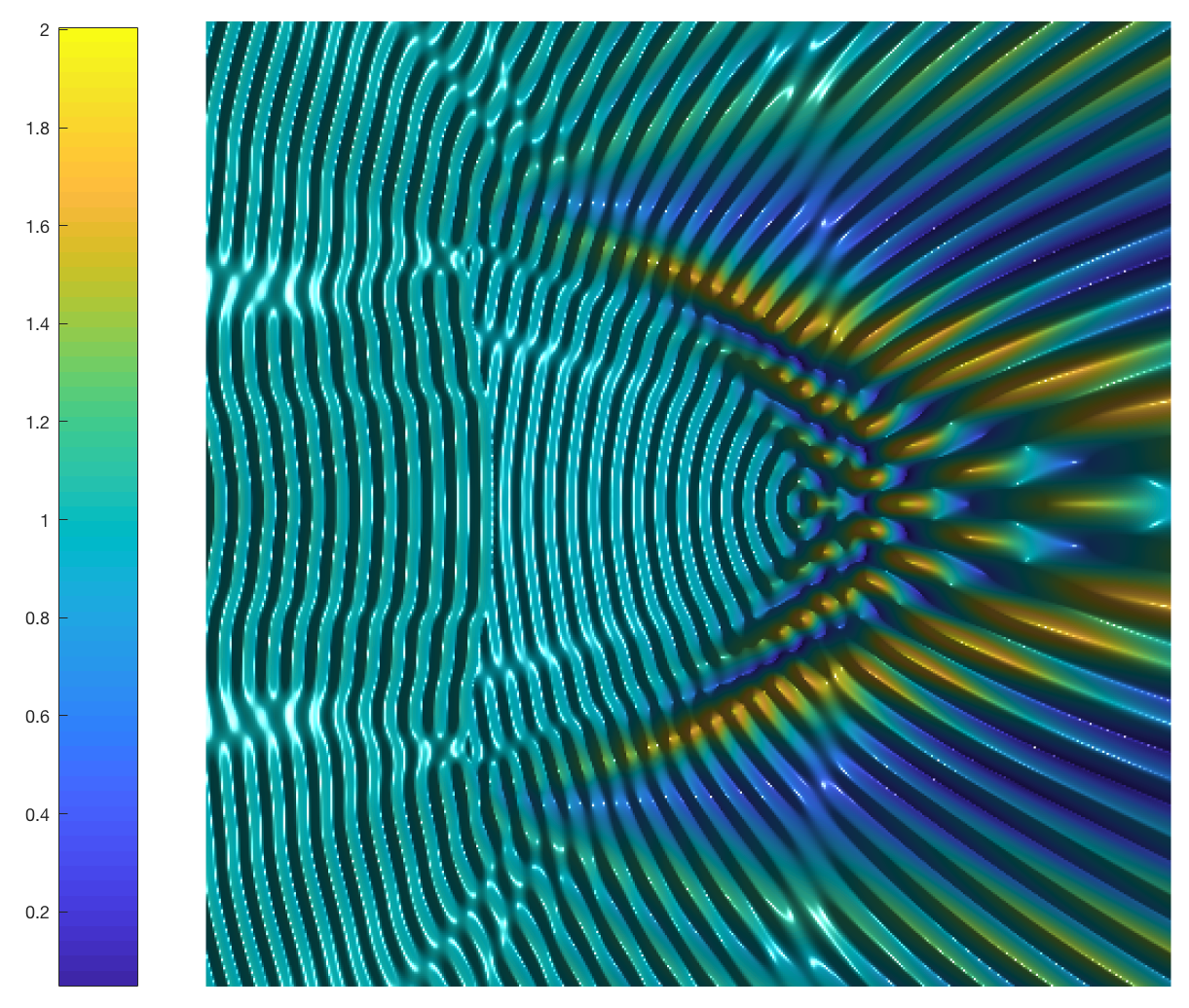}}
		\hfill
		\caption{Scattering of a plane wave $ \exp (i \kappa x)$ by a penetrable  bean shape scatterer of acoustical size $\kappa a =50$. We have taken computational grid of size $2\times{65} \times{129}+ 129 \times 65 $ and integral over base region is approximated by mens of Addition theorem  method.}
	\end{center}
\end{figure}

For a pictorial visualization, we plot the absolute part of the total and scattered fields in  Figure \ref{fig:-bean-K25}. In this experiment, we have taken $\kappa a =50$ and contrast function $m(\bsx)$ is $-0.3$ within inhomogeneity and zero elsewhere.  Again, for this experiment, the numerical solution is obtained using our approach with  Addition theorem where we have three digits of accuracy.


\begin{table} [!h]	
	\begin{center} 
		\begin{tabular}{c|c|c|c|c|c|c}\hline
			Grid Size & Unknowns  & \multicolumn{2}{c|}{$L^{\infty}$}&\multicolumn{2}{c|} {$L^2$} &  numIt\\ 
			\cline{3-7}&   & $\varepsilon_{\infty}$ & Order&  $\varepsilon_{2}$ & Order&  \\ 
			\hline
			$2\times9\times33+17\times17 $&379&9.88e-01&-&7.56e-01& - &13\\ 
			\hline
			$2\times17\times33 +33\times33$&1395&1.11e-01&3.16e+00&8.90e-02&3.09e+00&34\\ 
			\hline
			$2\times33\times33 +65\times65$&5347&5.63e-03&4.30e+00&4.52e-03&4.30e+00&51\\ 
			\hline
			$2\times65\times33 +129\times129$&20931&1.68e-04&5.07e+00&9.10e-05&5.63e+00&55\\ 
			\hline
			$2\times129\times33 +257\times257$&82819&3.70e-06&5.50e+00&1.63e-06&5.80e+00&58\\ 
		\end{tabular}  
		\caption{Convergence for the bean  scatterer with $ \kappa a= 10 \pi$ when pre-corrected trapezoidal rule is used for base integration.}
		\label{table:-PCTKA10PIBEAN}
	\end{center} 
\end{table}

\begin{table} [!h]	
	\begin{center} 
		\begin{tabular}{c|c|c|c|c|c|c}\hline
			Grid Size & Unknowns  & \multicolumn{2}{c|}{$L^{\infty}$}&\multicolumn{2}{c|} {$L^2$} &  numIt\\ 
			\cline{3-7}&   & $\varepsilon_{\infty}$ & Order&  $\varepsilon_{2}$ & Order&  \\ 
			\hline
			$2\times 9 \times 5+17\times 9 $&243&1.76e+00&-&1.42e+00& - &13\\ 
			\hline
			$2\times17 \times 9 +33\times 17$&867&5.56e-01&1.66e+00&3.72e-01&1.94e+00&26\\ 
			\hline
			$2\times33\times 17 +65 \times 33$&3267&4.51e-03&6.95e+00&4.70e-03&6.30e+00&37\\ 
			\hline
			$2\times65 \times33 +129\times 65$&12675&1.00e-04&5.49e+00&7.99e-05&5.88e+00&43\\ 
			\hline
			$2\times129\times 65 +257\times 129$&49923&1.62e-06&5.95e+00&1.25e-06&6.00e+00&48\\ 
		\end{tabular}  
	\caption{Convergence for the bean shape scatterer with $ \kappa a= 10 \pi$, when Addition theorem method is used for base integration.}
	\label{tabel:-ADTHKA10PIBEAN} 
	\end{center} 
\end{table} 
  
 \begin{exmp}(\textit{Computational efficiency})
 \end{exmp}	
In this example, we present numerical results to  corroborate the growth in computational complexity of our method.  As the thickness of the boundary region $\Omega_{B}$ is  kept to only a few wave lengths, therefore, transverse  integral in Eq. (\ref{eq:-BdrySplitK}) can be approximated accurately using a certain number of points independent of the wavelength. In view of this,  as discussed in section \ref{sec:-Comp-Cost}, the time complexity of our accelerated method exhibits $O(N \log N)$ growth. 

In order to demonstrate this, we again take $\Omega$ as a disc of unit radius with refractive index $n(\bsx)=\sqrt{2}$ if $\bsx \in \Omega$ and approximate the volume potential $\mathcal{A}(v)$ at all grid points in $\Omega_h$ at different levels of discretization. The computational results obtained by Addition theorem method for $\kappa a =30$, and those from pre-corrected trapezoidal rule for $\kappa a =50$, under the proposed framework, are reported in Table \ref{table:-ADTH_COST} and \ref{table:-PCT_COST} respectively. The numerical results clearly demonstrate that, for a fixed thin boundary region, we achieve
 the computational complexity of $O (N \log N)$.

\begin{table}[!h] 	
	\begin{center} 
		\begin{tabular}{c|c|c|c|c|c|c}\hline
			Grid Size & \multicolumn{2}{c|}{$L^{\infty}$}&\multicolumn{2}{c|} {$L^2$} & \multicolumn{2}{c} {Time (sec.)}  \\ 
			\cline{2-7}&    $\varepsilon_{\infty}$ & Order&  $\varepsilon_{2}$ & Order& accel & un-accel  \\ 			
			\hline
		
		$2\times33\times33 +65\times33$&7.24e-03&-&6.51e-03&-&5.00e+0&1.1e+01\\ 
		\hline
		$2\times65\times33 +129\times65$&1.10e-04&6.04e+00&1.36e-04&5.58e+00&1.70e+01&2.9e+01\\ 
		\hline
		$2\times129\times33 +257\times129$&4.22e-06&4.71e+00&4.97e-06&4.77e+00&4.20e+01&1.08e+02\\ 
		\hline
		$2\times257\times33 +513\times257$&1.26e-07&5.07e+00&1.18e-07&5.40e+00&1.72e+02&5.95e+02\\ 
			\hline 
		\end{tabular}  
		\caption{Performance of our method with Addition theorem  for a disc of size $\kappa a =30$ with refractive index $n(\bsx) =\sqrt{2}$. }
		\label{table:-ADTH_COST}
	\end{center} 
\end{table}

\begin{table}[!h] 	
	\begin{center} 
		\begin{tabular}{c|c|c|c|c|c|c}\hline
			Grid Size & \multicolumn{2}{c|}{$L^{\infty}$}&\multicolumn{2}{c|} {$L^2$} & \multicolumn{2}{c} {Time (sec.)}  \\ 
			\cline{2-7}&    $\varepsilon_{\infty}$ & Order&  $\varepsilon_{2}$ & Order& accel & un-accel  \\ 
		
		\hline
		$2\times33\times33 +65\times65$&8.84e-02&-&8.74e-02&-&7.0e+01&4.0e+00\\ 
		\hline
		$2\times65\times33 +129\times129$&6.81e-03&3.70e+00&4.42e-03&4.30e+00&1.50e+01&2.5e+01\\ 
		\hline
		$2\times129\times33 +257\times257$&2.27e-04&4.90e+00&1.54e-04&4.85e+00&3.90e+01&9.4e+01\\ 
		\hline
		$2\times257\times33 +513\times513$&3.41e-07&9.38e+00&2.62e-07&9.20e+00&1.54e+02&6.82e+02\\ 
			\hline 
		\end{tabular}  
		\caption{Performance of our method with pre-corrected trapezoidal rule  for a disc of size $\kappa a =50$ with refractive index $n(\bsx) =\sqrt{2}$. }
		\label{table:-PCT_COST}
	\end{center} 
\end{table}

%
%
%

 \begin{exmp}(\textit{Comparison with other high-order methods})
 \end{exmp}
In this example, we compare the performance of two algorithms implemented  under the present framework with existing high-order methods for discontinuous scattering media. To the best of our knowledge, high-order rates for discontinuous material interface  achieved by two solvers given in \cite{medvinsky2013high} and \cite{Anand2015highorder}.  In the aforementioned references, the first is based on the differential equation formulation, where as the second on the integral equation formulation. 
In addition, we have also included a comparison of our present method with those obtained by original Addition theorem method (ATM).

A comparative study of our  results with those of ATM and with the method
 in \cite{Anand2015highorder}, that we refer to as PUM, is given in Table \ref{table:-ADTHCOMP}.  In this table,  the results reported in first two rows are taken from Table $5.2$ on page $64$ in \cite{hyde2002fast}. Results reported in third and fourth rows are obtained by using our present approach with Addition theorem.   We see that, the present approach  produces substantially more accurate results than those obtained by ATM alone in \cite{hyde2002fast}. For instance, we  see in the second row that ATM produces an error $2.75\time 1.e-07$ using $2093000$ unknowns, whereas, our method produces an error $7.63\time 1.e-08$  while using only $79413$ unknowns.  In the rows five to eight, we have made similar comparisons with PUM. For this, the numerical results in  rows five and six  are taken form Table $1$ and $2$ at page $268$ in \cite{Anand2015highorder}, while, last two rows report numerical results obtained by our algorithm.  We observe that, the current method provides substantially more accurate results than those presented in \cite{Anand2015highorder}.  
   
	\begin{table}[!h]	
	\begin{center}
		\begin{tabular}{c|c|c|c|c|c} \hline
		Algorithm &	$\kappa a$  & $n$ & Unknown & $\varepsilon_{\infty}$ & $\varepsilon_{2}$ \\ 		
			\hline
			\hline
			ATM &4 $\pi$ &  $\sqrt{2}$  &$66000$ & 1.13e-03 & - \\
			ATM &4 $\pi$ &  $\sqrt{2}$  &$2093000$ & 2.75e-07 & - \\	
           Present &4 $\pi$ &  $\sqrt{2}$  &$12675$ & 1.70e-04& 7.33e-05 \\
			Present &4 $\pi$ &  $\sqrt{2}$  &$79413$ & 7.63e-08 & 1.14e-08 \\
			\hline	
		    PUM & 10 &  $\sqrt{3}+i\sqrt{2}$  &$25155$ & 3.94e-04 & 6.31e-04 \\	
			PUM & 10 &  $\sqrt{3}+i\sqrt{2}$  &$33411$ & 3.64e-05 & 5.66e-05 \\
			 Present& 10 &  $\sqrt{3}+i\sqrt{2}$  &$12675$ & 1.16e-04 & 8.09e-05 \\
			  Present& 10 &  $\sqrt{3}+i\sqrt{2}$  &$24303$ & 3.85e-06 & 2.63e-06 \\
			\hline
		\end{tabular}
		\caption{ Accuracy of the Addition theorem method when it used in the proposed framework.} 
		\label{table:-ADTHCOMP}
	\end{center}	
\end{table}

In \cite{medvinsky2013high}, a fourth-order algorithm for transmission scattering problem using method of difference potential is proposed. While this method converges to high-order even for discontinuous scatterers it does require a large number of unknowns to achieve reasonable accuracy.

To compare the performance of our method with the aforementioned algorithm, we compute scattering by an elliptical obstacle with major axis $b=1.8$ and minor axis $a=0.6$.  We take  $\kappa=10$ and $n(\bsx)=\sqrt{3}$ when $\bsx \in \Omega$ and one otherwise.  For this scattering problem,
in Table \ref{table:-MTDVC},  numerical results for the method of difference potential are taken from Table 6 in \cite{medvinsky2013high} and compared with those coming from our method. We see that the present algorithm produces better accuracy than those of  \cite{medvinsky2013high} while using smaller number of unknowns.  For example, we can see, in the last row of  Table \ref{table:-MTDVC}, our algorithm yields more than three times better accuracy while using $320$ times fewer unknowns.
\begin{table}[!h]	
	\begin{center}
		\begin{tabular}{c|c|c|c|c|c} \hline
			$\kappa $ & $n$ &\multicolumn{2}{c|}{Method of difference potential } & \multicolumn{2}{c}{current}  \\ 	
			
			\cline{3-6}& & Unknown & $\varepsilon_{\infty}$  & Unknown& $\varepsilon_{\infty}$\\ 			
			\hline
		
			10 &  $\sqrt{3}$ &$2097152$ & 3.35e-03  & 13715  &1.34e-03 \\		
		
			10 &$\sqrt{3}$ & $8388608$ & 2.09e-04 & 26194 &5.63 e-05  \\	
			\hline		
		\end{tabular}
		\caption{Comparison of our method with Method of Difference Potential proposed in \cite{medvinsky2013high}. }
		\label{table:-MTDVC}
	\end{center}	
\end{table}	

\section{Conclusions}
\label{sec:-Conclu}

In this paper, we have proposed an $O(N \log N)$ high-order method for the numerical solution of Lippmann-Schwinger integral equation in the two dimensions. More importantly, our methodology provides a framework that allows enhancement of convergence rates for those fast methods that converge rapidly for the case of smooth scattering media but yields poor convergence when the scattering medium has a jump discontinuity across the material interface, without adversely affecting their asymptotic computational cost. As illustrations, we have implemented two different $O(N\log N)$ algorithms under the proposed framework, namely, the Addition theorem method and the pre-corrected trapezoidal rule. We have shown an improved high-order convergence for both of these approaches.  Our numerical results clearly illustrate that, in case of discontinuous varying media, for given number of unknowns, accuracy obtained by our method compares 
well with those achieved by other existing solvers. 




\section*{Acknowledgments}

Akash Anand gratefully acknowledges support from IITK-ISRO Space Technology Cell through contract No.  STC/MATH/2014100.  Authors also thank Jagabandhu Paul for giving access to his implementation of the Addition Theorem based volumetric scattering solver.



\section*{References}

\bibliography{LS2DImproved.bib}
\bibliographystyle{elsarticle-num}

\end{document}